\documentclass[11pt,reqno]{amsart}
\usepackage{indentfirst, amssymb,amsmath,amsthm, mathrsfs,cite,graphicx,float}
\newtheorem*{theoA}{Theorem A}
\newtheorem*{theoB}{Theorem B}
\newtheorem*{theoC}{Theorem C}
\newtheorem*{theoD}{Theorem D}

\newtheorem{theo}{Theorem}[section]
\newtheorem{lem}{Lemma}[section]

\newtheorem{ex}{Example}[section]

\newtheorem*{Canonical function}{Canonical function}

\newcommand{\pa}{\partial}
\newcommand{\ol}{\overline}
\newcommand{\be}{\begin{equation}}
\newcommand{\ee}{\end{equation}}
\newcommand{\beas}{\begin{eqnarray*}}
\newcommand{\eeas}{\end{eqnarray*}}
\newcommand{\bea}{\begin{eqnarray}}
\newcommand{\eea}{\end{eqnarray}}
\newcommand{\bs}{\begin{small}}
\newcommand{\es}{\end{small}}
\renewcommand{\epsilon}{\varepsilon}
\numberwithin{equation}{section}

\begin{document}
\title[System of Fermat-type PDDE\lowercase{s}]
{Solutions of systems of certain Fermat-type PDDE\lowercase{s}}
\author[R. Biswas and R. Mandal]{Raju Biswas and Rajib Mandal}
\date{}
\address{Raju Biswas, Department of Mathematics, Raiganj University, Raiganj, West Bengal-733134, India.}
\email{rajubiswasjanu02@gmail.com}
\address{Rajib Mandal, Department of Mathematics, Raiganj University, Raiganj, West Bengal-733134, India.}
\email{rajibmathresearch@gmail.com}
\maketitle
\let\thefootnote\relax
\begin{abstract} The objective of this paper is to investigate the existence and the forms of the pair of finite order entire and meromorphic solutions of some certain systems of 
Fermat-type partial differential-difference equations of several complex variables. These results represent some refinements and generalizations of the earlier findings, especially the results due to Xu {\it et al.} (J. Math. Anal. Appl. 483(2) (2020)). We provide some examples to support the results.
\end{abstract}
\section{Introduction}
\noindent Nevanlinna's value distribution theory has occupied a central position in the field of complex analysis for more than a century. This theory is an important tool for studying the 
properties of meromorphic functions. The study of the Fermat-type functional equation has been a subject of interest in the field of complex analysis in the context of Nevanlinna's 
theory. The present-day central objective in the value distribution theory of Nevanlinna is to find the precise solutions
to various Fermat-type equations in complex variables. By a meromorphic function $f$ on $\mathbb{C}^n$ ($n\in\mathbb{N}$), we mean that $f$ can be written as a quotient of two holomorphic functions without common zero sets in $\mathbb{C}^n$.  Notationally, we write $f:=g/h$, where $g$ and $h$ are relatively prime holomorphic functions on $\mathbb{C}^n$ such that $h \not\equiv 0$ and $f^{-1}(\infty)\not=\mathbb{C}^n$. In particular, the entire function of several complex variables are holomorphic throughout $\mathbb{C}^n$. 
\subsection{Notations}
Let $z=(z_1,z_2,\ldots,z_n)\in\mathbb{C}^n$, $a\in\mathbb{C}\cup\{\infty\}$, $k\in\mathbb{N}$ and $r>0$. We consider some notations from \cite{14,S1,Y1}. Let $\ol B_n(r):=\{z\in\mathbb{C}^n: |z|\leq r\}$, where $|z|^2:=\sum_{j=1}^n|z_j|^2$. 
The exterior derivative splits $d:=\pa+\ol{\pa}$ and twists to $d^c:=\frac{i}{4\pi}(\ol\pa-\pa)$. The standard Kaehler metric on $\mathbb{C}^n$ is 
given by $v_n(z):=dd^c|z|^2$. Define $\omega_n(z):=dd^c\log |z|^2\geq 0$ and $\sigma_n(z):=d^c \log |z|^2\wedge \omega_n^{n-1}(z)$ on $\mathbb{C}^n\setminus\{0\}$. Thus $\sigma_n(z)$ defines a positive measure on $\pa B_n:=\{z\in\mathbb{C}^n: |z|= r\}$ with total measure $1$.
The zero-multiplicity of a holomorphic function $h$ at a point 
$z\in \mathbb{C}^n$ is defined to be the order of vanishing of $h$ at $z$ and denoted by $\mathcal{D}_h^0(z)$. A divisor of $f$ on $\mathbb{C}^n$ is an integer valued function which is locally the difference between the 
zero-multiplicity functions of $g$ and $h$ and it is denoted by $\mathcal{D}_f:=\mathcal{D}_g^0-\mathcal{D}_h^0$ (see, P. 381, \cite{51}). Let $a\in\mathbb{C}\cup\{\infty\}$ be such that $f^{-1}(a)\not=\mathbb{C}^n$. Then the $a$-divisor $\nu_f^a$ of $f$ 
is the divisor associated with the holomorphic functions $g-ah$ and $h$ (see, P. 346, \cite{14} and P. 12, \cite{12}). 
In \cite{Y1}, Ye has defined the counting function and the valence function with respect to $a$ respectively as follows:
\beas n(r,a,f):=r^{2-2n}\int_{S(r)} \nu_f^a v_n^{n-1}\;\text{and}\;N(r,a,f):=\int_0^r\frac{ n(r,a,f)}{t}dt. \eeas
We write 
\beas N(r,a,f)=\left\{\begin{array}{ll}
&N\left(r,\frac{1}{f-a}\right),\;\text{when}\;a\not=\infty\\[1mm]
&N(r,f),\;\text{when}\;a=\infty.\end{array}\right.\eeas 
The proximity function \cite{14,Y1} of $f$ is defined as follows :
\beas
&&m(r,f):=\int_{\pa B_n(r)} \log^+|f(z)|\sigma_n(z),\quad\text{when}\;a=\infty\\
&&m\left(r, \frac{1}{f-a}\right):=\int_{\pa B_n(r)} \log^+\frac{1}{|f(z)-a|}\sigma_n(z),\quad\text{when}\;a\not=\infty.\eeas
\noindent By denoting $S(r):= \ol{B}_n(r)\cap \text{supp}\;\nu_f^a$, where
$\text{supp}\;\nu_f^a=\ol{\left\{z\in\mathbb{C}^n: \nu_f^a(z)\not=0\right\}}$ (see, P. 346, \cite{14}). The notation $N_k\left(r,\frac{1}{f-a}\right)$ is known as truncated valence function. In particular, $N_1\left(r,\frac{1}{f-a}\right)=\ol N\left(r,\frac{1}{f-a}\right)$ is the truncated valence function of simple $a$-divisors of $f$ in $S(r)$. 
In $N_k\left(r,\frac{1}{f-a}\right)$, the $a$-divisors of $f$ in $S(r)$ of multiplicity $m$  are counted $m$-times if $m< k$ and $k$-times if $m\geq k$. The Nevanlinna characteristic function is defined by $T(r,f)=N(r,f)+m(r,f)$, which is increasing for $r$.
Given a meromorphic function $f$ on $\mathbb{C}^n$, we denote the order of $f$ by $\rho(f)$ such that 
\beas \rho(f)=\varlimsup\limits_{r\to \infty}\frac{\log^+ T(r,f)}{\log r},\;\;\text{where}\;\log^{+}x=\max\{\log x,0\}.\eeas
For further details, the reader is referred to \cite{51,12i,12,15,19,21,S1,Y1} and the references therein.\\[2mm]
\indent Given a meromorphic function $f(z)$ on $\mathbb{C}^n$, $f(z+c)$ is called a shift of $f$ and $\Delta (f)=f(z+c)-f(z)$ is called a difference operator of $f$, where 
$c(\not=0)\in\mathbb{C}^n$. An equation is called a partial differential equation (in brief, PDE) if the equation contains partial derivatives of f whereas if the equation also contains 
shifts or differences of $f$, then the equation is called a partial differential-difference equation (in brief, PDDE).
\subsection{Fermat-type equations and related results on $\mathbb{C}$}
\indent We now consider the Fermat-type equation
\bea\label{eq1.1} f^n(z)+g^n(z)=1,\quad\text{where}\quad n\in\mathbb{N}.\eea 
A significant number of researchers have demonstrated a keen interest in investigating
the Fermat-type equations (\ref{eq1.1}) for entire solutions \cite{49,39,7,17} and meromorphic solutions \cite{49,7,B1966} over
the past two decades by taking some variation of (\ref{eq1.1}). Yang and Li \cite{40} was the pioneer for introducing the study on transcendental meromorphic solutions of Fermat-type differential equation on $\mathbb{C}$. Liu \cite{Liu1} was the first who investigated on meromorphic solutions of Fermat-type difference equation as well as differential-difference equations on $\mathbb{C}$.\\[2mm]
In 2012, Liu et al. \cite{Liu2} investigated some Fermat-type differential-difference equations and obtained several results. Here we recall one of them.
\begin{theoA}\cite{Liu2} The finite order transcendental entire solutions of $\left(f'(z)\right)^2+f^2(z+c)=1$ must satisfy $f(z)=\sin(z\pm iB)$, where $c=2k\pi$ or $c=(2k+1)\pi$ with $k$ an integer and $B$ is a constant.\end{theoA}
We call the pair $(f,g)$ as a finite order meromorphic solutions for the system
\bea\label{14q1}\left\{\begin{array}{lll}
f^{m_1}+g^{n_1}=1,\\
f^{m_2}+g^{n_2}=1,\end{array}\right.\eea
if $f,g$ are meromorphic functions satisfying (\ref{14q1}) and $\rho(f,g)=\max\{\rho(f),\rho(g)\}<+\infty$.\\[1.5mm]
\indent In 2016, Gao \cite{Gao} initiated the study of a system of differential-difference equations and obtained the following result.
\begin{theoB}\cite{Gao} Suppose that $(f_1, f_2)$ is a pair of finite order transcendental entire solutions for the system of differential-difference equations
\beas\left\{\begin{array}{lll}
\left(f'_1(z)\right)^2+f_2^2(z+c)=1,\\
\left(f'_2(z)\right)^2+f_1^2(z+c)=1.\end{array}\right.\eeas
Then $(f_1, f_2)$ satisfies $(f_1,f_2)=\left(\sin(z-iA), \sin(z-iB)\right)$ or $(f_1,f_2)=\left(\sin(z+iA),\right.\\\left. \sin(z+iB)\right)$, where $c=k\pi$ with $k$ an integer and $A,B$ are constants.
\end{theoB}
Shortly afterwards, Gao \cite{Gao1} investigated the existence and the forms of the finite order entire solutions of the system of differential-difference equations
\bea\label{gao}\left\{\begin{array}{lll}
\left(f'_1(z)\right)^{m_1}+f_2^{n_1}(z+c)=Q_1(z),\\
\left(f'_2(z)\right)^{m_2}+f_1^{n_2}(z+c)=Q_2(z),\end{array}\right.\eea
where $m_i, n_i\in\mathbb{N}$ $(i=1,2)$, $Q_1(z)$ and $Q_2(z)$ are non-zeros polynomials and obtained the following results.
\begin{theoC}\cite{Gao1} The system (\ref{gao}) does not have any pair of finite order transcendental entire solutions $(f_1,f_2)$ whenever 
$n_1n_2>m_1m_2$ or $n_j>m_j/(m_j-1)$ for $m_j\geq 2$, $j=1,2$.\end{theoC}
\begin{theoD}\cite{Gao1} Let $(f_1,f_2)$ be the finite order transcendental entire solutions of the system (\ref{gao}) with $m_i=n_i=2$ $(i=1,2)$.  Then $Q_1(z)= A_1B_1$, $Q_2(z)=A_2B_2$, and
\beas f_1(z)=\left(A_1e^{az+b_1}-B_1e^{-az-b_1}\right)/(2a)\;f_2(z)=\left(A_2e^{az+b_2}-B_2e^{-az-b_2}\right)/(2a),\eeas
where $a^4=1$, $b_1$, $b_2$, $A_i(\not=0)$, $B_i(\not=0)$ $(i=1,2)$ are constants.\end{theoD}
\subsection{Fermat-type equations in several complex variables}
The basic conclusions about the solutions of the Fermat-type equation (\ref{eq1.1}) on $\mathbb{C}$ were also extended to the case of several complex variables (see \cite[\textrm{Theorem 2.3}]{700} \cite[\textrm{Theorem 1.3}]{20}).
Researchers have recently directed their efforts toward investigating the Fermat-type PDEs for entire and meromorphic solutions. Let \bea\label{esl}\sum_{i=1}^n\left(\frac{\pa u}{\pa z_i}\right)^m=1\eea
be the certain non-linear first order PDE introducing from the analogy with the Fermat-type equation $\sum_{i=1}^n\left(f_i\right)^m=1$, where $u:\mathbb{C}^n\to \mathbb{C}$, $z_i\in\mathbb{C}$, $f_i:\mathbb{C}\to \mathbb{C}$, and $m,n\geq 2$.
In 1999, Saleeby \cite{27} first started to study about the solutions of the Fermat-type PDEs and obtained the results for entire solutions of (\ref{esl}) on $\mathbb{C}^2$. Afterwards, Li \cite{28} extended these results to $\mathbb{C}^n$. In 2008, Li \cite{16} considered the equation (\ref{eq1.1}) with $n=2$ and showed that meromorphic solutions $f$ and $g$ of that equation on 
$\mathbb{C}^2$ must be constant if and only if $\pa f/\pa z_2$ and $\pa g/\pa z_1$ have the same zeros (counting multiplicities).
If $f=\pa u/\pa z_1$ and $g=\pa u/\pa z_2$, then any entire solutions of the PDE $f^2+g^2=1$ on $\mathbb{C}^2$ are necessarily linear \cite{13}.\\[1.5mm]
\indent In 2018, Xu and Cao \cite{23,26} first considered both difference operators and differential operators in Fermat-type equations of two complex variables and obtained the results about the existence and forms of transcendental entire solutions of finite order.
In 2020, Xu {\it et al.} \cite{200} significantly changed the research of this direction to Fermat-type system of PDDEs. Actually, the authors \cite{200} considered the following systems of PDDEs
\bea&&\label{14q}\left\{\begin{array}{lll}
\left(\frac{\pa f_1(z_1,z_2)}{\pa z_1}\right)^{n_1}+f_2(z_1+c_1,z_2+c_2)^{m_1}=1,\\
\left(\frac{\pa f_2(z_1,z_2)}{\pa z_1}\right)^{n_2}+f_1(z_1+c_1,z_2+c_2)^{m_2}=1,\end{array}\right.\\\text{and}
&&\label{tre}\left\{\begin{array}{lll}
\left(\frac{\pa f_1(z_1,z_2)}{\pa z_1}\right)^2+\left(f_2(z_1+c_1,z_2+c_2)-f_1(z_1,z_2)\right)^2=1,\\
\left(\frac{\pa f_2(z_1,z_2)}{\pa z_1}\right)^2+\left(f_1(z_1+c_1,z_2+c_2)-f_2(z_1,z_2)\right)^2=1,\end{array}\right.
\eea
where $(c_1,c_2)\in\mathbb{C}^2$ and $m_j,n_j\in\mathbb{N}$ for $j=1,2$ and obtained the existence and explicit representations of transcendental entire solutions with finite order for systems (\ref{14q}) and (\ref{tre}) separately.\\[2mm]
\indent For the recent developments on the solutions of Fermat-type PDDEs, we refer to \cite{32,MB1,35,Hal1,Raj1,500,Xu1, BM2024,601,7,49,28,16,Liu2,602,603,604,600,RR,XLX2024} and the references therein.
\subsection{ Basic Notations}
Let $I=(i_1,i_2,\ldots,i_k)\in\mathbb{Z}^k_+$ be a multi-index with length $\Vert I\Vert=\sum_{j=1}^k i_j$ and 
\beas \pa^I f=\frac{\pa^{\Vert I\Vert} f}{\pa z_1^{i_1}\cdots\pa z_k^{i_k}}.\eeas 
It is evident that $\mathcal{P}(z)=\sum_{\Vert I\Vert=0}^\alpha \beta_I z_1^{i_1}\cdots z_n^{i_n}$ is a polynomial of degree $\alpha$ in several complex variables, where 
$\beta_{I}\in\mathbb{C}$ such that $\beta_{I}$ are not all zero at a time for $\Vert I\Vert=\alpha$. Suppose that $\mathcal{P}(z+c)-\mathcal{P}(z)\equiv B\in\mathbb{C}$, for 
any $c\in\mathbb{C}^n\setminus\{(0,0,\ldots,0)\}$ and 
$\mathcal{P}(z)=\sum_{j=1}^n a_jz_j+\mathcal{Q}(z)+A$, where $A\in\mathbb{C}$, $\deg (\mathcal{Q}(z))\geq 2$. Now, $\mathcal{P}(z+c)-\mathcal{P}(z)\equiv B$ implies 
that $\sum_{j=1}^n a_j c_j+\mathcal{Q}(z+c)-\mathcal{Q}(z)\equiv B$. For $c\in\mathbb{C}^n\setminus\{(0,0,\ldots,0)\}$, it is evident that for any monomials $M(z)$ of several 
complex variables with $\deg(M(z))\geq 2$, we have $M(z+c)- M(z)\not \equiv \text{constant}$. Thus, we have $\mathcal{Q}(z+c)\equiv \mathcal{Q}(z)$ and $\sum_{j=1}^n a_jc_j=B$.\\[2mm]
\indent Now we can rewrite the polynomial $\mathcal{Q}(z)$ with $\deg(\mathcal{Q}(z))\geq 2$ in such a way that it has the terms from the polynomials like 
$\Psi(\epsilon_{j_1}z_{j_1}+\epsilon_{j_2}z_{j_2}+\ldots+\epsilon_{j_m} z_{j_m})$ in $\epsilon_{j_1}z_{j_1}+\epsilon_{j_2}z_{j_2}+\ldots+\epsilon_{j_m} z_{j_m}$ such that $\epsilon_{j_1}c_{j_1}+\epsilon_{j_2}c_{j_2}+\ldots+\epsilon_{j_m} c_{j_m}=0$, $\epsilon_{j_1}$, $\cdots$, $\epsilon_{j_m}\in\mathbb{C}$ ($1\leq j_1,j_2,\cdots,j_m\leq n$) and $\deg(\Psi)\geq 2$. Since $\mathcal{Q}(z)$ is periodic, so we can express $\mathcal{Q}(z)$ as
\bea\label{K1} \mathcal{Q}(z)=\sum_{\lambda}Q_\lambda(z)\quad\text{and}\quad Q_\lambda(z)=\prod_{\alpha} Q_\alpha (z),\eea
where $\lambda$ belongs to the finite index set $I_1$ of the family $\{Q_\lambda(z) : \lambda \in I_1\}$ and $\alpha$ belongs to the finite index set $I_2$ of the family 
$\{Q_\alpha(z) : \alpha \in I_2\}$ with
\beas Q_\alpha(z)&=&\sum_{\substack{j_1,j_2=1,\\j_1<j_2}}^n \Psi_{2,\alpha,j_1,j_2}(t_{j_1}z_{j_1}+t_{j_2}z_{j_2})+\sum_{\substack{j_1,j_2,j_3=1,\\j_1<j_2<j_3}}^n \Psi_{3,\alpha,j_1,j_2,j_3}(\zeta_{j_1}z_{j_1}+\zeta_{j_2}z_{j_2}+\zeta_{j_3} z_{j_3})\nonumber\\
&&+\cdots+\sum_{\substack{j_1,j_2,\ldots,j_n=1,\\j_1<j_2<\ldots<j_n}}^n \Psi_{n,\alpha,j_1,j_2,\ldots,j_n}(\eta_{j_1}z_{j_1}+\eta_{j_2}z_{j_2}+\cdots+\eta_{j_n} z_{j_n})\nonumber\eeas
where $\eta_i,\zeta_i,t_i, A\in\mathbb{C}$ $(1\leq i\leq n)$, $\deg \mathcal{Q}(z)=\deg \mathcal{P}(z)$, $\Phi_{m,\alpha,j_1,j_2,\ldots,j_m}(\eta_{j_1}z_{j_1}+\eta_{j_2}z_{j_2}+\ldots+\eta_{j_m} z_{j_m})$ is a polynomial in $\eta_{j_1}z_{j_1}+\eta_{j_2}z_{j_2}+\ldots+\eta_{j_m} z_{j_m}$. Here $\eta_i,\zeta_i,t_i\in\mathbb{C}$ $(1\leq i\leq n)$ are chosen from the conditions $t_{j_1}c_{j_1}+t_{j_2}c_{j_2}=0$, $\zeta_{j_1}c_{j_1}+\zeta_{j_2}c_{j_2}+\zeta_{j_3} c_{j_3}=0$, $\eta_{j_1}c_{j_1}+\eta_{j_2}c_{j_2}+\ldots+\eta_{j_m} c_{j_m}=0$. It is also clear that for $ j_1=1,2,\ldots,n$, we have
\beas\frac{\pa Q_\alpha(z)}{\pa z_{j_1}}&=&t_{j_1}\sum_{\substack{j_1,j_2=1,\\j_1<j_2}}^n \Psi_{2,\alpha,j_1,j_2}'(t_{j_1}z_{j_1}+t_{j_2}z_{j_2})\\
&&+\zeta_{j_1}\sum_{\substack{j_1,j_2,j_3=1,\\j_1<j_2<j_3}}^n \Psi_{3,\alpha,j_1,j_2,j_3}'(\zeta_{j_1}z_{j_1}+\zeta_{j_2}z_{j_2}+\zeta_{j_3} z_{j_3})\\
&&+\ldots+\eta_{j_1}\sum_{\substack{j_1,j_2,\ldots,j_n=1,\\j_1<j_2<\ldots<j_n}}^n \Psi_{n,\alpha,j_1,j_2,\ldots,j_n}'(\eta_{j_1}z_{j_1}+\eta_{j_2}z_{j_2}+\cdots+\eta_{j_n} z_{j_n}).\hspace{5cm}\eeas
Let $F_{f}(z)$ denote the partial differential function of a finite order transcendental meromorphic function $f$ with $N(r,f)=S(r,f)$ involving $n(\in\mathbb{N})$ different homogeneous terms such that
\bea\label{R1} F_{f}(z)=\sum_{m=1}^{n}\sum_{\Vert I\Vert=m}a_{I}(z)\pa^I f(z)=\sum_{m=1}^{n}\sum_{\Vert I\Vert=m}a_{I}(z)\frac{\pa^{\Vert I\Vert} f}{\pa z_1^{i_1}\cdots\pa z_k^{i_k}},\eea
where $z=\left(z_1,z_2,\ldots,z_n\right)$ and $a_{I}(z)$ are small functions of $f(z)$ on several complex variables such that $a_I(z)$ are not all identically zero at a time.\\[1.5mm]
\indent 
Motivated by the results of \cite{32,200, Gao,Gao1} and due to continue the research for further investigations on equations of several complex variables, in this paper, we 
consider the following systems of Fermat-type PDDEs of several complex variables for $1\leq \mu\leq n$:
\bea&&\label{fg}\left\{\begin{array}{lll} \left(F_{f_1}(z)\right)^{m_1}+P_1(z)f_2^{n_1}(z+c)=Q_1(z),\\[2mm]
\left(F_{f_2}(z)\right)^{m_2}+P_2(z)f_1^{n_2}(z+c)=Q_2(z),\end{array}\right.\\[1mm]
&&\label{e1}\left\{\begin{array}{lll}
\left(a_1\frac{\pa f_1(z)}{\pa z_\mu}\right)^2+\left(a_2f_1(z)+a_3f_2(z+c)+a_4\frac{\pa^2 f_1(z)}{\pa z_\mu^2}\right)^2=1,\\[2mm]
\left(a_1\frac{\pa f_2(z)}{\pa z_\mu}\right)^2+\left(a_2f_2(z)+a_3f_1(z+c)+a_4\frac{\pa^2 f_2(z)}{\pa z_\mu^2}\right)^2=1,\end{array}\right.\\[1mm]\text{and}
&&\label{e4}\left\{\begin{array}{lll}
\left(a_1\frac{\pa f_1(z)}{\pa z_1}+a_2\frac{\pa f_1(z)}{\pa z_2}+\cdots+a_n\frac{\pa f_1(z)}{\pa z_n}\right)^2+\left(a_{n+1}f_1(z)+a_{n+2}f_2(z+c)\right)^2=1,\\[2mm]
\left(a_1\frac{\pa f_2(z)}{\pa z_1}+a_2\frac{\pa f_2(z)}{\pa z_2}+\cdots+a_n\frac{\pa f_2(z)}{\pa z_n}\right)^2+\left(a_{n+1}f_2(z)+a_{n+2}f_1(z+c)\right)^2=1,\end{array}\right.\eea
where $m_j,n_j\in\mathbb{N}$ $(j=1,2)$, $c=\left(c_1,c_2,\ldots,c_n\right)\in\mathbb{C}^n$, $a_j\in\mathbb{C}\setminus\{0\}$ $(1\leq j\leq n+2)$ and $P_i(z)(\not\equiv0)$, $Q_i(z)(\not\equiv 0)$ are small functions of $f_1$ and $f_2$ on $\mathbb{C}^n$ for $i=1,2$. Throughout the paper, we denotes
\bs\bea\label{K2}\left\{\begin{array}{lll}
 y=\left(a_\mu z_1-a_1z_\mu, \ldots, a_\mu z_{\mu-1}-a_{\mu-1}z_\mu, a_\mu z_{\mu+1}-a_{\mu+1}z_\mu,\ldots,a_\mu z_n-a_n z_\mu\right),\\[2mm] 
 s=\left(a_\mu c_1-a_1c_\mu, \ldots, a_\mu c_{\mu-1}-a_{\mu-1}c_\mu, a_\mu c_{\mu+1}-a_{\mu+1}c_\mu,\ldots,a_\mu c_n-a_n c_\mu\right),\\[1.5mm]
y_1=(z_1,z_2,\ldots,z_{\mu-1},z_{\mu+1},\ldots,z_n), s_1=(c_1,c_2,\ldots,c_{\mu-1},c_{\mu+1},\ldots,c_n)\\[1.5mm]
\gamma_1(k)=(a_{k+1}c_\mu-ia_1)/(2a_1) \quad\text{and}\quad \gamma_2(k)=(a_{k+1}c_\mu+ia_1)/(2a_1).\end{array}\right.\eea\es
In our all statements, we assume that $c=(c_1,c_2,\ldots,c_n)\in\mathbb{C}^n\setminus\{(0,0,\ldots,0)\}$. To the best of our knowledge, the above system of equations has not been considered before. 
\section{Main Results}
In the following result, we investigate the existence of solutions of the system (\ref{fg}).
\begin{theo}\label{th} The system (\ref{fg}) does not have any pair of finite order transcendental meromorphic solutions $(f_1,f_2)$ on $\mathbb{C}^n$ with $N(r,f_j)=S(r,f_j)$ whenever 
$n_1n_2>m_1m_2$ or $n_j>m_j/(m_j-1)$ for $m_j\geq 2$, $j=1,2$.
\end{theo}
\noindent Note that, in a particular case, \textrm{Theorem \ref{th}} becomes the results of \cite{32,200} and thus \textrm{Theorem \ref{th}} is the generalization of the results \cite{32,200}. \\[1mm]
\noindent For the pair of finite order transcendental entire solutions of the system (\ref{e1}), we obtain the following results.
\begin{theo}\label{T1} If $a_2=\pm a_3$, then $$f_1(z)=z_\mu/\sqrt{a_1^2+a_2^2c_\mu^2}+h_1(y_1),\quad f_2(z)=z_\mu/\sqrt{a_1^2+a_2^2c_\mu^2}+h_2(y_1)$$
are finite order transcendental entire solutions of the system (\ref{e1}), where $h_j(y_1)$ $(j=1,2)$ are finite order transcendental entire functions with periods $2s_1$.\end{theo}
\noindent To make it easier to write down the following \textrm{Theorems \ref{T11}-\ref{T13}}, we  introduced such a condition ($\mathcal{A}$):\\[1.5mm]
 ($\mathcal{A}$) The  constants $b_j,A,B,K_i\in\mathbb{C}$ $(1\leq j\leq n\quad\text{and}\quad 1\leq i\leq 4)$ such that $K_1K_2=1=K_3K_4$, $\Psi_1(z)$ is a polynomial 
defined in (\ref{K1}) with $\Psi_1(z)\equiv 0$, if $\Psi_1(z)$ contain the variable $z_\mu$ and $g_j(y_1)$ $(1\leq j\leq 6)$ are finite order entire functions satisfying 
\beas &&a_2g_1(y_1)+a_3g_2(y_1+s_1)=\gamma_1(1)K_1e^{\sum_{j=1,j\not=\mu}^nb_jz_j+A}+\gamma_2(1)K_2e^{-\sum_{j=1,j\not=\mu}^nb_jz_j-A},\\[1mm]
&&a_2g_2(y_1)+a_3g_1(y_1+s_1)=\gamma_1(1)K_3e^{-\sum_{j=1,j\not=\mu}^nb_jz_j+B}+\gamma_2(1)K_4e^{\sum_{j=1,j\not=\mu}^nb_jz_j-B},\\[1mm]
&&a_2g_3(y_1)+a_3g_4(y_1+s_1)\equiv0\quad\text{and}\quad a_2g_4(y_1)+a_3g_3(y_1+s_1)\equiv0,\eeas
where $\gamma_1(1)$ and $\gamma_2(1)$ are defined in (\ref{K2}).
\begin{theo}\label{T11} If $a_2^2=\pm a_3^2$, then 
\beas f_1(z)&=&\frac{z_\mu}{2a_1}\left(K_1e^{\sum_{j=1,j\not=\mu}^nb_jz_j+\Psi_1(z)+A}+K_2e^{-\sum_{j=1,j\not=\mu}^nb_jz_j-\Psi_1(z)-A}\right)+g_1(y_1),\\
f_2(z)&=&\frac{z_\mu}{2a_1}\left(K_3e^{-\sum_{j=1,j\not=\mu}^nb_jz_j-\Psi_1(z)+B}+K_4e^{\sum_{j=1,j\not=\mu}^nb_jz_j+\Psi_1(z)-B}\right)+g_2(y_1)\eeas
are finite order transcendental entire solutions of the system (\ref{e1}) for which the condition ($\mathcal{A}$) holds and 
\beas\left\{\begin{array}{lll}
e^{2\sum_{j=1,j\not=\mu}^nb_jc_j}\equiv \left(a_3/a_2\right)^2\equiv e^{-2\sum_{j=1,j\not=\mu}^nb_jc_j},\quad e^{2A+2B}\equiv a_3^2K_2K_4/(a_2^2K_1K_3),\\[1.5mm]
e^{-\sum_{j=1,j\not=\mu}^nb_jc_j+A+B}\equiv -a_3K_4/(a_2K_1), \quad e^{\sum_{j=1,j\not=\mu}^nb_jc_j-A-B}\equiv -a_3K_3/(a_2K_2),\\[1.5mm]
e^{\sum_{j=1,j\not=\mu}^nb_jc_j+A+B}\equiv -a_3K_2/(a_2K_3), \quad e^{-\sum_{j=1,j\not=\mu}^nb_jc_j-A-B}\equiv -a_3K_1/(a_2K_4).\end{array}\right.\eeas.\end{theo} 
\begin{theo}\label{T12} If $a_1^4a_2^2=a_3^4a_4^2$, then 
\beas f_1(z)&=&\frac{K_1}{2a_1b_\mu}e^{\sum_{j=1}^nb_jz_j+\Psi_1(z)+A}-\frac{K_2}{2a_1b\mu}e^{-\sum_{j=1}^nb_jz_j-\Psi_1(z)-A}+g_3(y_1),\\
f_2(z)&=&-\frac{K_3}{2a_1b_\mu}e^{-\sum_{j=1}^nb_jz_j-\Psi_1(z)+B}+\frac{K_4}{2a_1b_\mu}e^{\sum_{j=1}^nb_jz_j+\Psi_1(z)-B}+g_4(y_1)\eeas
are finite order transcendental entire solutions of the system (\ref{e1}) for which the condition ($\mathcal{A}$) holds and 
\beas\left\{\begin{array}{lll}
e^{2\sum_{j=1}^nb_jc_j}\equiv \left(a_3/(a_1b_\mu)\right)^2\equiv e^{-2\sum_{j=1}^nb_jc_j}, &e^{2A+2B}\equiv a_3^2K_2K_4/(a_1^2b_\mu^2K_1K_3)\\[1.5mm]
e^{-\sum_{j=1}^nb_jc_j+A+B}\equiv ia_3K_4/(a_1b_\mu K_1), &e^{\sum_{j=1}^nb_jc_j-A-B}\equiv ia_3K_3/(a_1b_\mu K_2),\\[1.5mm]
e^{\sum_{j=1}^nb_jc_j+A+B}\equiv -ia_3K_2/(a_1b_\mu K_3), &e^{-\sum_{j=1}^nb_jc_j-A-B}\equiv -ia_3K_1/(a_1b_\mu K_4).\end{array}\right.\eeas\end{theo} 
\begin{theo}\label{T13} If $ia_1b_\mu+a_4b_\mu^2+a_2=(-1)^\nu a_3$ $(\nu=1,2)$, then
\beas f_1(z)&=&\frac{K_1}{2a_1b_\mu}e^{\sum_{j=1}^nb_jz_j+\Psi_1(z)+A}-\frac{K_2}{2a_1b_\mu}e^{-\sum_{j=1}^nb_jz_j-\Psi_1(z)-A}+g_3(y_1),\\
f_2(z)&=&\frac{K_3}{2a_1b_\mu}e^{\sum_{j=1}^nb_jz_j+\Psi_1(z)+B}-\frac{K_4}{2a_1b_\mu}e^{-\sum_{j=1}^nb_jz_j-\Psi_1(z)-B}+g_4(y_1)\eeas
are finite order transcendental entire solutions of the system (\ref{e1}) for which the condition ($\mathcal{A}$) holds and 
\beas\left\{\begin{array}{llll}
e^{2\sum_{j=1}^nb_jc_j}\equiv 1, \quad e^{2A-2B}\equiv K_2K_3/(K_1K_4), &e^{\sum_{j=1}^nb_jc_j-A+B}\equiv (-1)^{\nu+1}K_4/K_2,\\[1.5mm]
 e^{\sum_{j=1}^nb_jc_j+A-B}\equiv (-1)^{\nu+1}K_2/K_4, &e^{-\sum_{j=1}^nb_jc_j+A-B}\equiv (-1)^{\nu+1} K_3/K_1, \\[1.5mm]
 e^{-\sum_{j=1}^nb_jc_j-A+B}\equiv (-1)^{\nu+1} K_1/K_3.\end{array}\right.\eeas \end{theo} 
The following example related to \textrm{Theorem \ref{T1}} is reasonable.
\begin{ex} Let 
$f_1(z)=\frac{1}{6}e^{3z_1-3z_2+3z_3}+\frac{1}{6}e^{-3z_1+3z_2-3z_3}+e^{-\frac{2\log (-6)}{\pi i}(z_2+z_3)}$ and $f_2(z)=\frac{i}{6}e^{-3z_1+3z_2-3z_3}-\frac{i}{6}e^{3z_1-3z_2+3z_3}-e^{-\frac{2\log (-6)}{\pi i}(z_2+z_3)}$. Then $\rho(f_1,f_2)=1$ and $(f_1,f_2)$ satisfies the system (\ref{e1}) with $a_1=i$, $a_2=-18$, $a_3=3$, $a_4=2$ and $c=(\pi i, -\pi i, \pi i/2)$.
\end{ex}
For the pair of finite order transcendental entire solutions of the system (\ref{e4}), we obtain the following result.
\begin{theo}\label{T2} If $a_{n+1}=\pm a_{n+2}$, then $f_1(z)=z_\mu/\sqrt{a_1^2+a_2^{n+1}c_\mu^2}+g_1(y)$, $f_2(z)=z_\mu/\sqrt{a_1^2+a_{n+1}^2c_\mu^2}+g_2(y)$
are finite order transcendental entire solutions of the system (\ref{e4}), where $g_j(y)$ $(j=1,2)$ are finite order transcendental entire functions with periods $2s$ satisfying 
$\sum_{k=1}^n a_k\frac{\pa g_j(y)}{\pa z_k}\equiv 0$.\end{theo}
\noindent To make it easier to write down the following  \textrm{Theorems \ref{T22}-\ref{T24}}, we  introduced such a condition ($\mathcal{B}$):\\[1.5mm]
($\mathcal{B}$) The  constants $b_j,A,B,K_i\in\mathbb{C}$ $(1\leq j\leq n\quad\text{and}\quad 1\leq i\leq 4)$ such that $K_1K_2=1=K_3K_4$ and $h_j(y)$ $(1\leq j\leq 6)$ are finite order entire functions satisfying 
\beas&& a_{n+1}h_1(y)+a_{n+2}h_2(y+s)\equiv\gamma_1(n)K_1e^{\sum_{j=1}^nb_jz_j+A}+\gamma_2(n)K_2e^{-\sum_{j=1}^nb_jz_j-A},\\[1mm]
&& a_{n+1}h_2(y)+a_{n+2}h_1(y+s)\equiv\gamma_1(n)K_3e^{-\sum_{j=1}^nb_jz_j+B}+\gamma_2(n)K_4e^{\sum_{j=1}^nb_jz_j-B},\\[1mm]
&&a_{n+1}h_3(y)+a_{n+2}h_4(y+s)\equiv 0,\quad a_{n+1}h_4(y)+a_{n+2}h_3(y+s)\equiv 0\\[1mm]
&& a_{n+1}h_5(y)+a_{n+2}h_6(y+s)\equiv\gamma_1(n)K_1e^{\sum_{j=1}^nb_jz_j+A}+\gamma_2(n)K_2e^{-\sum_{j=1}^nb_jz_j-A}\\[1mm]\text{and}
&& a_{n+1}h_6(y)+a_{n+2}h_5(y+s)\equiv\gamma_1(n)K_3e^{\sum_{j=1}^nb_jz_j+B}+\gamma_2(n)K_4e^{-\sum_{j=1}^nb_jz_j-B},\eeas
where $\gamma_1(n)$ and $\gamma_2(n)$ are defined in (\ref{K2}).
\begin{theo}\label{T22} If $a_{n+1}^2=\pm a_{n+2}^2$, then 
\beas f_1(z)&=&\frac{z_\mu}{2a_\mu}\left(K_1e^{\sum_{j=1}^nb_jz_j+A}+K_2e^{-\sum_{j=1}^nb_jz_j-A}\right)+h_1\left(y\right),\\ 
f_2(z)&=&\frac{z_\mu}{2a_\mu}\left(K_3e^{-\sum_{j=1}^nb_jz_j+B}+K_4e^{\sum_{j=1}^nb_jz_j-B}\right)+h_2(y)\eeas
are finite order transcendental entire solutions of the system (\ref{e4}) for which the condition ($\mathcal{B}$) holds and 
\beas\left\{\begin{array}{lll}
e^{2\sum_{j=1}^nb_jc_j}\equiv \left(a_{n+2}/a_{n+1}\right)^2\equiv e^{-2\sum_{j=1}^nb_jc_j}, \quad  e^{2A+2B }\equiv a_{n+2}^2K_2K_4/(a_{n+1}^2K_1K_3),\\
e^{-\sum_{j=1}^nb_jc_j+A+B}\equiv -a_{n+2}K_4/(a_{n+1}K_1), \quad  e^{\sum_{j=1}^nb_jc_j+A+B}\equiv -a_{n+2}K_2/(a_{n+1}K_3),\\
e^{\sum_{j=1}^nb_jc_j-A-B}\equiv -a_{n+2}K_3/(a_{n+1}K_2), \quad e^{-\sum_{j=1}^nb_jc_j-A-B}\equiv -a_{n+2}K_1/(a_{n+1}K_4).\end{array}\right.\eeas\end{theo} 
\begin{theo}\label{T23} If $a_{n+1}\not=(-1)^\nu a_{n+2}$ $(\nu=1,2)$, then 
\beas f_1(z)&=&\frac{K_1e^{\sum_{j=1}^nb_jz_j+A}-K_2e^{-\sum_{j=1}^nb_jz_j-A}}{2\sum_{j=1}^n a_j b_j}+h_3(y),\\
f_2(z)&=&\frac{K_3e^{\sum_{j=1}^nb_jz_j+B}-K_4e^{-\sum_{j=1}^nb_jz_j-B}}{2\sum_{j=1}^n a_j b_j}+h_4(y)\eeas
are finite order transcendental entire solutions of the system (\ref{e4}) for which the conditions ($\mathcal{B}$) and $\sum_{j=1}^n a_jb_j=i(a_{n+1}-(-1)^\nu a_{n+2})$ hold with 
\bea\label{J3}\left\{\begin{array}{lll}
e^{2\sum_{j=1}^nb_jc_j}\equiv 1, e^{2A-2B}\equiv K_2K_3/(K_1K_4), e^{\sum_{j=1}^nb_jc_j-A+B}\equiv (-1)^{\nu+1}K_4/K_2,\\[1.5mm]
e^{\sum_{j=1}^nb_jc_j+A-B}\equiv (-1)^{\nu+1}K_2/K_4, e^{-\sum_{j=1}^nb_jc_j+A-B}\equiv (-1)^{\nu+1}K_3/K_1,\\[1.5mm]
e^{-\sum_{j=1}^nb_jc_j-A+B}\equiv (-1)^{\nu+1}K_1/K_3.\end{array}\right.\eea\end{theo} 
\begin{theo}\label{T24} If $a_{n+1}=(-1)^\nu a_{n+2}$ $(\nu=1,2)$, then 
\beas f_1(z)&=&\frac{z_\mu}{2a_\mu}\left(K_1e^{\sum_{j=1}^nb_jz_j+A}+K_2e^{-\sum_{j=1}^nb_jz_j-A}\right)+h_5(y),\\
f_2(z)&=&\frac{z_\mu}{2a_\mu}\left(K_3e^{\sum_{j=1}^nb_jz_j+B}+K_4e^{-\sum_{j=1}^nb_jz_j-B}\right)+h_6(y)\eeas
are finite order transcendental entire solutions of the system (\ref{e4}) for which the conditions ($\mathcal{B}$) and $\sum_{j=1}^n a_jb_j=0$ hold with (\ref{J3}).
\end{theo} 
The following example related to \textrm{Theorem \ref{T2}} is reasonable.
\begin{ex} Let 
\bs\beas f_1(z)&=&-\frac{1}{60}e^{3z_1+z_2-2z_3+5z_4}+\frac{1}{60}e^{-3z_1-z_2+2z_3-5z_4}+e^{\frac{15\log (2/3)}{8\pi i}(-8z_1+z_2+z_3+z_4)},\\ 
f_2(z)&=&\frac{i}{60}e^{3z_1+z_2-2z_3+5z_4}-\frac{i}{60}e^{-3z_1-z_2+2z_3-5z_4}-e^{\frac{15\log (2/3)}{8\pi i}(-8z_1+z_2+z_3+z_4)}.\eeas\es
Then $\rho(f_1,f_2)=1$ and $(f_1,f_2)$ satisfies the system of PDDEs for $c=(\pi i/3, 2\pi i, \pi i, \pi i/5)$, 
\bs\beas\left\{\begin{array}{lll}
\left(\frac{\pa f_1(z)}{\pa z_1}-2\frac{\pa f_1(z)}{\pa z_2}+3\frac{\pa f_1(z)}{\pa z_3}+7\frac{\pa f_1(z)}{\pa z_4}\right)^2+\left(-12i f_1(z)-18 i f_2(z+c)\right)^2=1,\\
\left(\frac{\pa f_2(z)}{\pa z_1}-2\frac{\pa f_2(z)}{\pa z_2}+3\frac{\pa f_2(z)}{\pa z_3}+7\frac{\pa f_2(z)}{\pa z_4}\right)^2+\left(-12i f_2(z)-18i f_1(z+c)\right)^2=1.\end{array}\right.\eeas\es
\end{ex}
\noindent In the proof of the main results of the paper, we have used the Nevanlinna’s theory of several complex variables, the difference analogue of the lemma on the logarithmic derivative in several complex variables \cite{4,14} and the Lagrange's auxiliary equations \cite[Chapter 2]{430} for quasi-linear PDEs.
\section {Key Lemmas} 
The following lemmas are essential to prove the main results of this paper.
\begin{lem}\label{lem1}\cite[\textrm{Lemma 1.5}, P. 239]{12i}\cite[\textrm{Lemma 3.1, P. 211}]{12} Let $f_j\not\equiv 0\;(j=1, 2,3)$ be meromorphic functions on $\mathbb{C}^n$ such that $f_1$ is not constant and $f_1+f_2+f_3\equiv 1$ with
\bs\beas \sum_{j=1}^3\left\{N_2(r, 0; f_j)+2\ol N(r, f_j)\right\}<\lambda T(r, f_1)+O(\log^+T(r, f_1))\eeas\es
holds as $r\to\infty$ out side of a possible exceptional set of finite linear measure,
where $\lambda<1$ is a positive number. Then either $f_2\equiv 1$ or $f_3\equiv 1$. \end{lem}
\begin{lem}\label{lem2}\cite[Theorem 4.3.4, P. 247]{19} Let $f(z)$ be an entire function on $\mathbb{C}^n$ such that $f(0)\not=0$ and $\rho(N(r,0,f))<\infty$. 
Then there exists an entire function $g(z)$ and a canonical function $\phi(z)$ (see \cite[Theorem 4.3.2, P. 245]{19})
such that $f(z)=\phi(z)e^{g(z)}$.
\end{lem}
\begin{lem}\label{lem8}\cite[\textrm{Lemma 2.1}, P. 282]{12ai} \cite[\textrm{Lemma 3.58}]{12} If $g$ is a transcendental entire function on $\mathbb{C}^n$ and if $f$ is a meromorphic function of positive order on $\mathbb{C}$, then $f\circ g$ is of infinite order.\end{lem}
\begin{lem}\label{lem9}\cite[\textrm{Proposition 3.2}, P.240]{12ii} \cite[\textrm{Lemma 3.59}]{12} Let $P$ be a non-constant entire function on $\mathbb{C}^n$. Then
\beas \rho(e^P)=\left\{\begin{array}{clcr}\deg(P)&:\text{if}\; P \;\text{is a polynomial},\\
+\infty&:\text{otherwise}.\end{array}\right.\eeas
\end{lem}
\begin{lem}\label{lem4}\cite[\textrm{Theorem 2.1}, P. 242]{12i}\cite[\textrm{Lemma 1.106}]{12}  Suppose that $a_0(z), a_1(z),\ldots$,\\$a_m(z)$ $(m\geq 1)$ are meromorphic functions on $\mathbb{C}^n$ and $g_0(z), g_1(z),
\ldots,g_m(z)$ are entire functions on $\mathbb{C}^n$ such that $g_j(z)-g_k(z)$ are not constants for $0\leq j < k \leq n$. If $\sum_{j=0}^na_j(z)e^{g_j(z)}\equiv 0$ and $T(r, a_j)=o(T(r))$, $j=0,1,\ldots,n$ hold as $r\to\infty$ out side of a possible exceptional set of finite linear measure, where $T(r)=\min\limits_{0\leq j < k \leq n} T (r, e^{g_j-g_k})$,
then $a_j(z) \equiv 0$ $(j =0,1,2,\ldots,n)$.\end{lem}
\begin{lem}\label{lem20}\cite[\textrm{Lemma 3.2}, P. 385]{51} Let $f$ be a non-constant meromorphic function on $\mathbb{C}^n$. Then for any $I\in \mathbb{Z}_+^n$, $T(r, \pa^I f)=O(T(r,f))$ 
for all $r$ except possibly a set of finite Lebesgue measure, where $I=\left(i_1,i_2,\ldots,i_n\right)\in\mathbb{Z}_+^n$ denotes a multiple index with $\Vert I\Vert=i_1+i_2+\cdots+i_n$, $\mathbb{Z}_+=\mathbb{N}\cup\{0\}$.  
\end{lem}
\begin{lem}\label{lem5}\cite{4,14} Let $f$ be a non-constant meromorphic function with finite order on $\mathbb{C}^n$ such that $f(0)\not=0,\infty$. Then for $c\in\mathbb{C}^n$, $m\left(r, \frac{f(z)}{f(z+c)}\right)+m\left(r, \frac{f(z+c)}{f(z)}\right)=S(r, f)$
holds for all $r>0$ outside of a possible exceptional set $E\subset [1,\infty)$ of finite logarithmic measure $\int_E \frac{dt}{t}<\infty$. \end{lem}
\begin{lem}\label{lem6}\cite{3,25} Let $f$ be a non-constant meromorphic function with finite order on $\mathbb{C}^n$ and $I=(i_1,i_2,\ldots,i_n)$ be a multi-index with length $\Vert I\Vert=\sum_{j=1}^n i_j$. Assume that $T(r_0, f)\geq e$ for some $r_0$. Then $m\left(r, \frac{\pa^I f}{f}\right)=S(r,f)$, holds for all $r\geq r_0$ outside a set $E\subset (0,\infty)$ of finite logarithmic measure $\int_E dt/t<\infty$.\end{lem}
\section{Proofs of the main theorems}  
\begin{proof}[\bf{Proof of Theorem \ref{th}}] 
Let $(f_1,f_2)$  be a pair of finite order transcendental meromorphic functions satisfies (\ref{fg}) with $N(r, f_j)=S(r, f_j)$ for $j=1,2$. 
Then $F_{f_j}(z), f_j(z+c)$ ($j=1,2$) are finite order transcendental meromorphic functions on several complex variables with $N(r,F_{f_j}(z))=S(r,f_j)=N(r,f_j(z+c))$, where $F_{f}(z)$ is defined in (\ref{R1}). 
Now we discuss the following two cases.\\
{\bf Case 1.}  When $n_1n_2>m_1m_2$, then by using \textrm{Lemma \ref{lem5}}, we have 
\bea \label{fg1}T(r, f_j(z))=m(r, f_j(z))+S(r,f_j)
&\leq&m\left(r, \frac{f_j(z)}{ f_j(z+c)}\right)+m(r,  f_j(z+c))+S(r,f_j)\nonumber\\
 &\leq& T(r, f_j(z+c))+S(r,f_j),\eea
 holds for all $r>0$ outside of a possible exceptional set $E=E_1\cup E_2\subset [1,\infty)$ of finite logarithmic measure. In view of Valiron-Mokhon'ko lemma [page. 29, \cite{12}], \textrm{Lemma \ref{lem6}} and (\ref{fg1}), we deduce that
\bea\label{14p1} n_1T(r, f_2(z))&\leq& n_1T(r, f_2(z+c))+S(r,f_2)=T\left(r, f_2^{n_1}(z+c)\right)+S(r,f_2)\nonumber\\[1mm]
&\leq& T\left(r, P_1(z)f_2^{n_1}(z+c)\right)+S(r,f_1)+S(r,f_2)\nonumber\\[1mm]
&=&T\left(r, \left(F_{f_1}(z)\right)^{m_1}-Q_1(z) \right)+S(r,f_1)+S(r,f_2)\nonumber\\[1mm]
&\leq&m_1T\left(r,F_{f_1}(z)\right)+S(r,f_1)+S(r,f_2)\nonumber\\[1mm]
&=&m_1m\left(r,F_{f_1}(z)\right)+S(r,f_1)+S(r,f_2)\nonumber\\
&\leq&m_1m\left(r,\sum\limits_{m=1}^{n}\sum\limits_{\Vert I\Vert=m}a_{I}(z)\frac{\pa^I f_1(z)}{f_1(z)}\right)+m_1m(r, f_1(z))\nonumber\\
&&+S(r,f_1)+S(r,f_2)\nonumber\\[1mm]
&\leq&m_1m(r, f_1(z))+S(r,f_1)+S(r,f_2)\nonumber\\[1mm]
&\leq& m_1T(r, f_1(z))+S(r,f_1)+S(r,f_2),\eea
holds for all $r>0$ outside of a possible exceptional set $E$ of finite logarithmic measure. Similarly, we have 
\bea\label{14p2}  n_2T(r, f_1(z))\leq m_2T(r, f_2(z))+S(r,f_1)+S(r,f_2), \;r\not\in E.\eea
From (\ref{14p1}) and (\ref{14p2}), we have 
\beas \left(n_1n_2-m_1m_2\right)T(r, f_j(z))\leq S(r,f_1)+S(r,f_2), \;r\not\in E,\eeas
which arise a contradiction, since $f_1,f_2$ are finite order transcendental meromorphic functions and $n_1n_2>m_1m_2$.\\
{\bf Case 2.} When $n_j>m_j/(m_j-1)$ for $m_j\geq 2$, $j=1,2$.
By Nevanlinna second fundamental theorem for small functions, \textrm{Lemma \ref{lem5}} and (\ref{fg}), we have
\bea\label{fg2} m_1T(r, F_{f_1}(z))&=&T\left(r, \left(F_{f_1}(z)\right)^{m_1}\right)+S(r,f_1)\nonumber\\[1mm]
&\leq&\ol N\left(r, \left(F_{f_1}(z)\right)^{m_1}\right)+\ol N\left(r, 0;\left(F_{f_1}(z)\right)^{m_1}\right)\nonumber\\[1mm]
&&+\ol N\left(r,0;\left(F_{f_1}(z)\right)^{m_1}-Q_1(z)\right)+S(r,f_1)\nonumber\\[1mm]
&=&\ol N\left(r, 0;F_{f_1}(z)\right)+\ol N\left(r,0; P_1(z)f_2^{n_1}(z+c)\right)+S(r,f_1)\nonumber\\[1mm]
&\leq& T(r,F_{f_1}(z))+\ol N\left(r, 0; f_2(z+c)\right)+S(r,f_1)+S(r,f_2)\nonumber\\[1mm]\text{\it i.e.,}\quad 
(m_1-1)T(r, F_{f_1}(z))&\leq& T(r, f_2(z+c))+S(r,f_1)+S(r,f_2).\eea
Similarly, we have
\bea\label{14p3}(m_2-1)T(r, F_{f_2}(z))&\leq& T(r, f_1(z+c))+S(r,f_1)+S(r,f_2).\eea  
Again, in view of Valiron-Mokhon'ko lemma [page. 29, \cite{12}], \textrm{Lemma \ref{lem5}}, (\ref{fg}), (\ref{fg2}) and (\ref{14p3}), we have
\bea\label{14p4} n_1T(r, f_2(z+c))&=& T(r, P_1(z)f_2^{n_1}(z+c))+S(r,f_1)+S(r,f_2)\nonumber\\[1mm]
&=&T\left(r, \left(F_{f_1}(z)\right)^{m_1}-Q_1(z)\right)+S(r,f_2)\nonumber\\[1mm]
&\leq&m_1T(r, F_{f_1}(z))+S(r,f_1)+S(r,f_2)\nonumber\\[1mm]
&\leq&\frac{m_1}{m_1-1}T(r, f_2(z+c))+S(r, f_1)+S(r,f_2)\nonumber\\[1mm]\text{\it i.e.,}\quad
\left(n_1-\frac{m_1}{m_1-1}\right)T(r, f_2(z+c))&\leq&S(r, f_1)+S(r,f_2).\eea 
Similarly, we have
\bea\label{14p5} \left(n_2-\frac{m_2}{m_2-1}\right)T(r, f_1(z+c))&\leq&S(r, f_1)+S(r,f_2).\eea 
From (\ref{14p4}) and (\ref{14p5}), we arise at a contradiction, since $f_1,f_2$ are finite order transcendental meromorphic functions and $n_j>m_j/(m_j-1)$ for $m_j\geq 2$, $j=1,2$. This completes the proof. 
\end{proof}
\begin{proof}[\bf{Proof of Theorem \ref{T1}}] {\bf Part 1.} The first part is common for all Theorems \ref{T1}-\ref{T13}.
Let $(f_1,f_2)$ be a pair of finite order transcendental entire functions  satisfies the system (\ref{e1}).
Now (\ref{e1}) can be written as 
\beas\left\{\begin{array}{lll}\prod\limits_{j=1}^2\left(a_1\frac{\pa f_1(z)}{\pa z_\mu}-(-1)^ji\left( a_2f_1(z)+a_3f_2(z+c)+a_4\frac{\pa^2 f_1(z)}{\pa z_\mu^2}\right)\right)=1\\[1mm]
\prod\limits_{j=1}^2\left(a_1\frac{\pa f_2(z)}{\pa z_\mu}-(-1)^ji\left( a_2f_2(z)+a_3f_1(z+c)+a_4\frac{\pa^2 f_2(z)}{\pa z_\mu^2}\right)\right)=1.\end{array}\right. \eeas
Here $a_1\frac{\pa f_j(z)}{\pa z_\mu}\pm i\left( a_2f_j(z)+a_3f_k(z+c)+a_4\frac{\pa^2 f_j(z)}{\pa z_\mu^2}\right)$ are finite order transcendental entire functions for $1\leq j,k\leq 2$ with $ j\not=k$ and have no zeros on $\mathbb{C}^n$. In view of the \textrm{Lemma \ref{lem2}}, we have
\bea\label{14q2}\left\{\begin{array}{llll}
&a_1\frac{\pa f_1(z)}{\pa z_\mu}+i\left( a_2f_1(z)+a_3f_2(z+c)+a_4\frac{\pa^2 f_1(z)}{\pa z_\mu^2}\right)=K_1e^{P(z)},\\[2mm]
&a_1\frac{\pa f_1(z)}{\pa z_\mu}-i\left( a_2f_1(z)+a_3f_2(z+c)+a_4\frac{\pa^2 f_1(z)}{\pa z_\mu^2}\right)=K_2e^{-P(z)},\\[2mm]
&a_1\frac{\pa f_2(z)}{\pa z_\mu}+i\left( a_2f_2(z)+a_3f_1(z+c)+a_4\frac{\pa^2 f_2(z)}{\pa z_\mu^2}\right)=K_3e^{Q(z)},\\[2mm]
&a_1\frac{\pa f_2(z)}{\pa z_\mu}-i\left( a_2f_2(z)+a_3f_1(z+c)+a_4\frac{\pa^2 f_2(z)}{\pa z_\mu^2}\right)=K_4e^{-Q(z)},\end{array}\right.\eea
where $K_1,K_2,K_3,K_4\in\mathbb{C}\setminus\{0\}$ such that $K_1K_2=1=K_3K_4$ and $P(z),Q(z)$ are entire functions on $\mathbb{C}^n$. From (\ref{14q2}), we deduce that 
\bea\label{eq1}\left\{\begin{array}{llll}
&a_1\frac{\pa f_1(z)}{\pa z_\mu}=\frac{K_1e^{P(z)}+K_2e^{-P(z)}}{2},\\[1mm]
&a_2f_1(z)+a_3f_2(z+c)+a_4\frac{\pa^2 f_1(z)}{\pa z_\mu^2}=\frac{K_1e^{P(z)}-K_2e^{-P(z)}}{2i},\\[1mm]
&a_1\frac{\pa f_2(z)}{\pa z_\mu}=\frac{K_3e^{Q(z)}+K_4e^{-Q(z)}}{2},\\[1mm]
&a_2f_2(z)+a_3f_1(z+c)+a_4\frac{\pa^2 f_2(z)}{\pa z_\mu^2}=\frac{K_3e^{Q(z)}-K_4e^{-Q(z)}}{2i}.\end{array}\right.\eea
Since $\rho(f_i(z))<\infty$ $(i=1,2)$, in view of \textrm{Lemmas \ref{lem8}}, \ref{lem9} and \ref{lem20}, we get from (\ref{eq1}) that $P(z)$ and $Q(z)$ are polynomials on $\mathbb{C}^n$. 
The following cases arise separately in the proofs of all \textrm{Theorems} \ref{T1}-\ref{T13}.\\[1.5mm]
{\bf Part 2.} Now we start to prove \textrm{Theorem \ref{T1}} properly.
Let $P(z),Q(z)$ be both constants. From (\ref{eq1}), we have
\bea\label{er4}\left\{\begin{array}{lll} 
&a_1\frac{\pa f_1(z)}{\pa z_\mu}=A_1,& a_2f_1(z)+a_3f_2(z+c)+a_4\frac{\pa^2 f_1(z)}{\pa z_\mu^2}=A_2,\\[2mm]
&a_1\frac{\pa f_2(z)}{\pa z_\mu}=A_3,& a_2f_2(z)+a_3f_1(z+c)+a_4\frac{\pa^2 f_2(z)}{\pa z_\mu^2}=A_4,\end{array}\right.\eea 
where $A_j\in\mathbb{C}$ for $1\leq j\leq 4$ with $A_1^2+A_2^2=1$ and $A_3^2+A_4^2=1$.
From (\ref{er4}), we have 
\beas f_1(z)=(A_1/a_1)z_\mu+h_1(y_1)\quad\text{and}\quad f_2(z)=(A_3/a_1)z_\mu+h_2(y_1),\eeas
where $h_j(y_1)$ $(j=1,2)$ are finite order transcendental entire functions of $z_1,z_2,\ldots,z_{\mu-1}$, $z_{\mu+1},\ldots,z_n$. Thus, we deduce that 
\beas
&&((a_2A_1+a_3A_3)/a_1)z_\mu+(a_2h_1(y_1)+a_3h_2(y_1+s_1))+a_3c_\mu A_3/a_1\equiv A_2\\[1mm]\text{and}
&&((a_2A_3+a_3A_1)/a_1)z_\mu+(a_2h_2(y_1)+a_3h_1(y_1+s_1))+a_3c_\mu A_1/a_1\equiv A_4.
\eeas
Since $h_j(y_1)$ $(j=1,2)$ are finite order transcendental entire functions, so we have 
\beas &&a_2A_1+a_3A_3= 0, a_2 A_3+a_3 A_1= 0, a_2h_1(y_1)+a_3h_2(y_1+s_1)\equiv 0\\[1mm]
&&a_2h_2(y_1)+a_3h_1(y_1+s_1)\equiv 0,a_3c_\mu A_3= a_1A_2\quad\text{and}\quad a_3c_\mu A_1= a_1 A_4.\eeas
For non-zero solution of system $a_2 A_1+a_3A_3=0$, $a_2A_3+a_3A_1=0$, we must have 
$\begin{vmatrix} 
a_2&a_3\\
a_3&a_2
\end{vmatrix}=0$, {\it i.e.,} $a_2=\pm a_3$, which implies that $A_1=\pm A_3$. It is easy to see that $A_1/A_2=-a_1/(a_2c_\mu)$. From 
$A_1^2+A_{2}^2=1$, we have $A_2=\pm a_2c_\mu/\sqrt{a_1^2+a_2^2c_\mu^2}$, $A_1=\mp a_1/\sqrt{a_1^2+a_2^2c_\mu^2}$ and $A_3=a_1/\sqrt{a_1^2+a_2^2c_\mu^2}$. Therefore, we have 
\beas f_1(z)=\frac{ z_\mu}{\sqrt{a_1^2+a_2^2c_\mu^2}}+h_1(y_1)\quad\text{and}\quad f_2(z)=\frac{z_\mu}{\sqrt{a_1^2+a_2^2c_\mu^2}}+h_2(y_1),\eeas
where $h_j(y_1)$ $(j=1,2)$ are finite order transcendental entire functions with periods $2s_1$.\end{proof}
\begin{proof}[\bf{Proof of Theorem \ref{T11}}] Let either $P(z)$ or $Q(z)$ be a constant. Let us assume that $P(z)$ be a constant and $Q(z)$ be a non-constant polynomial. From (\ref{eq1}), we have
$a_1\frac{\pa f_1(z)}{\pa z_\mu}=A_1$, $a_2f_1(z)+a_3f_2(z+c)+a_4\frac{\pa^2 f_1(z)}{\pa z_\mu^2}=A_2$, where $A_1,A_2\in\mathbb{C}$ with $A_1^2+A_2^2=1$. 
Now, $a_2f_1(z)+a_3f_2(z+c)=A_2$ implies that 
\beas \frac{\pa f_2(z+c)}{\pa z_\mu}=-\frac{a_2}{a_3}\frac{\pa f_1(z)}{\pa z_\mu}=-\frac{a_2}{a_3}\frac{A_1}{a_1},\eeas 
which contradicts that $\pa f_2(z)/\pa z_\mu$ is a transcendental entire function.\\[1.5mm]
Let $P(z),Q(z)$ be both non-constant polynomials.
By differentiating partially with respect to $z_\mu$ on both sides of the first equation in (\ref{eq1}), we have
\bea\label{eq3}\frac{\pa^2 f_1(z)}{\pa z_\mu^2}=\frac{K_1e^{P(z)}-K_2e^{-P(z)}}{2a_1}\frac{\pa P(z)}{\pa z_\mu}.\eea  
With the help of (\ref{eq3}), we deduce from the second equation of (\ref{eq1}) that
\bea\label{eq4}
&&a_2f_1(z)+a_3f_2(z+c)=\frac{K_1e^{P(z)}-K_2e^{-P(z)}}{2}\left(\frac{1}{i}-\frac{a_4}{a_1}\frac{\pa P(z)}{\pa z_\mu}\right).\eea
Differentiating partially (\ref{eq4}) with respect to $z_\mu$, we get
\bea\label{eq5}&&a_2\frac{\pa f_1(z)}{\pa z_\mu}+a_3\frac{\pa f_2(z+c)}{\pa z_\mu}=K_1e^{P(z)}\left(\frac{1}{2i}\frac{\pa P(z)}{\pa z_\mu}-\frac{a_4}{2a_1}\left(\frac{\pa P(z)}{\pa z_\mu}\right)^2-\frac{a_4}{2a_1}\frac{\pa^2P(z)}{\pa z_\mu^2}\right)\nonumber\\
&&+K_2e^{-P(z)}\left(\frac{1}{2i}\frac{\pa P(z)}{\pa z_\mu}-\frac{a_4}{2a_1}\left(\frac{\pa P(z)}{\pa z_\mu}\right)^2+\frac{a_4}{2a_1}\frac{\pa^2P(z)}{\pa z_\mu^2}\right). \eea
From (\ref{eq1}) and (\ref{eq5}), we obtain 
\bea
\label{tr}&&\Gamma_1(z)e^{P(z)+Q(z+c)}+\Psi_1(z) e^{-P(z)+Q(z+c)}-\frac{K_3}{K_4}e^{2Q(z+c)}\equiv 1, \eea
where
\bea\label{yt} &&\Gamma_1(z)=\frac{a_1K_1}{a_3K_4}\left(\frac{1}{i}\frac{\pa P(z)}{\pa z_\mu}-\frac{a_4}{a_1}\left(\frac{\pa P(z)}{\pa z_\mu}\right)^2-\frac{a_4}{a_1}\frac{\pa^2P(z)}{\pa z_\mu^2}-\frac{a_2}{a_1}\right)\\\text{and}
\label{yt1}&&\Psi_1(z)=\frac{a_1K_2}{a_3K_4}\left(\frac{1}{i}\frac{\pa P(z)}{\pa z_\mu}-\frac{a_4}{a_1}\left(\frac{\pa P(z)}{\pa z_\mu}\right)^2+\frac{a_4}{a_1}\frac{\pa^2P(z)}{\pa z_\mu^2}-\frac{a_2}{a_1}\right).\eea
In view of the first, third and fourth equations of (\ref{eq1}) and by means of similar argument as above, we have
\bea\label{14p10}&& \Gamma_2(z)e^{Q(z)+P(z+c)}+\Psi_2(z) e^{-Q(z)+P(z+c)}-\frac{K_1}{K_2}e^{2P(z+c)}\equiv 1,\\\text{where} 
\label{14yt} &&\Gamma_2(z)=\frac{a_1K_3}{a_3K_2}\left(\frac{1}{i}\frac{\pa Q(z)}{\pa z_\mu}-\frac{a_4}{a_1}\left(\frac{\pa Q(z)}{\pa z_\mu}\right)^2-\frac{a_4}{a_1}\frac{\pa^2Q(z)}{\pa z_\mu^2}-\frac{a_2}{a_1}\right)\\\text{and}
\label{14yt1}&&\Psi_2(z)=\frac{a_1K_4}{a_3K_2}\left(\frac{1}{i}\frac{\pa Q(z)}{\pa z_\mu}-\frac{a_4}{a_1}\left(\frac{\pa Q(z)}{\pa z_\mu}\right)^2+\frac{a_4}{a_1}\frac{\pa^2Q(z)}{\pa z_\mu^2}-\frac{a_2}{a_1}\right).\eea
From (\ref{tr}), it is clear that both $\Gamma_1(z)$ and $\Psi_1(z)$ are not simultaneously identically zero, otherwise we get a contradiction. Let $\Gamma_1(z)\equiv 0$ and $\Psi_1(z)\not\equiv 0$. From (\ref{tr}), we have  
\be\label{14p11} \Psi_1(z) e^{-P(z)+Q(z+c)}-\frac{K_3}{K_4}e^{2Q(z+c)}\equiv 1,\;\;\text{\it{i.e.,}}\;\;\Psi_1(z) e^{Q(z+c)}-\frac{K_3}{K_4}e^{2Q(z+c)+P(z)}-e^{P(z)}\equiv 0.\ee
From (\ref{14p11}), it is evident that $Q(z+c)-P(z)$ is a non-constant polynomial, otherwise a contradiction arises.
We claim that $Q(z+c)+P(z)$ and $2Q(z+c)+P(z)$ are non-constant polynomials. If possible, let $Q(z+c)+P(z)\equiv t_1$ and $2Q(z+c)+P(z)\equiv t_2$, where $t_1,t_2\in\mathbb{C}$. For above two situations, we deduce from (\ref{14p11}) that 
\bea\label{14p13}\begin{array}{llll} 
&\Psi_1(z) e^{t_1}-\frac{K_3}{K_4}e^{2t_1}\equiv e^{2P(z)},
&\Psi_1(z) e^{2Q(z+c)-t_2}-\frac{K_3}{K_4}e^{Q(z+c)}-e^{-Q(z+c)}\equiv 0.\end{array}\eea
In both circumstances, we get a contradiction from (\ref{14p13}) by using \textrm{Lemma \ref{lem4}}. Hence $Q(z+c)+P(z)$ and $2Q(z+c)+P(z)$ are non-constant polynomials.
In view of \textrm{Lemma \ref{lem4}}, we get a contradiction from (\ref{14p11}). Using similar argument, we again get contradictions from (\ref{tr}) and (\ref{14p10}), whenever $\Gamma_1(z)\not\equiv 0$, $\Psi_1(z)\equiv 0$; $\Gamma_2(z)\equiv 0$, $\Psi_2(z)\not\equiv 0$ and  $\Gamma_2(z)\not\equiv 0$, $\Psi_2(z)\equiv 0$. From (\ref{tr}), it is easy to see that
\bs\beas &&N\left(r, \Gamma_1(z)e^{P(z)+Q(z+c)}\right)=N\left(r, 0; \Gamma_1(z)e^{P(z)+Q(z+c)}\right)=N\left(r, \Psi_1(z)e^{-P(z)+Q(z+c)}\right)\\
&&=N\left(r, 0; \Psi_1(z)e^{-P(z)+Q(z+c)}\right)=N\left(r, -K_3e^{2Q(z+c)}/K_4\right)=N\left(r, 0;-K_3e^{2Q(z+c)}/K_4\right)\\
&&=S\left(r,-K_3e^{2Q(z+c)}/K_4\right).\eeas\es
In view of \textrm{Lemma \ref{lem1}} and (\ref{tr}), we have
\beas\text{either}\quad\Gamma_1(z)e^{P(z)+Q(z+c)}\equiv 1\quad\text{ or}\quad \Psi_1(z) e^{-P(z)+Q(z+c)}\equiv1,\eeas
where $\Gamma_1(z)$ and $\Psi_1(z)$ are given in (\ref{yt}) and (\ref{yt1}) respectively. It is evident from (\ref{14p10}) that
\bs\beas &&N\left(r, \Gamma_2(z)e^{Q(z)+P(z+c)}\right)=N\left(r, 0; \Gamma_2(z)e^{Q(z)+P(z+c)}\right)=N\left(r, \Psi_2(z)e^{-Q(z)+P(z+c)}\right)\\
&&=N\left(r, 0; \Psi_2(z)e^{-Q(z)+P(z+c)}\right)=N\left(r, -K_1e^{2P(z+c)}/K_2\right)=N\left(r, 0;-K_1e^{2P(z+c)}/K_2\right)\\
&&=S\left(r,-K_1e^{2P(z+c)}/K_2\right).\eeas\es
By applying \textrm{Lemma \ref{lem1}} in (\ref{14p10}), we get 
\beas\text{either}\quad\Gamma_2(z)e^{Q(z)+P(z+c)}\equiv 1\quad\text{or}\quad\Psi_2(z) e^{-Q(z)+P(z+c)}\equiv 1,\eeas
where $\Gamma_2(z)$ and $\Psi_2(z)$ are given in (\ref{14yt}) and (\ref{14yt1}) respectively.
Now the following cases arise.\\[2mm]
Let \bea\label{bv1}\begin{array}{lll}
\Gamma_1(z)e^{P(z)+Q(z+c)}\equiv 1\;\;\text{and}\;\;\Gamma_2(z)e^{Q(z)+P(z+c)}\equiv 1.\end{array}\eea
Using (\ref{bv1}), we deduce from (\ref{tr}) and (\ref{14p10}) respectively
\bea\label{bv5}\begin{array}{lll}
K_4\Psi_1(z) e^{-P(z)-Q(z+c)}/K_3\equiv1\quad\text{and}\quad K_2\Psi_2(z) e^{-Q(z)-P(z+c)}/K_1\equiv 1.\end{array}\eea
From (\ref{bv1}), it is clear that $P(z)+Q(z+c)$ and $Q(z)+P(z+c)$ are both constants, say $\chi_1$ and $\chi_2$ respectively, where $\chi_1,\chi_2\in\mathbb{C}$. Now 
$P(z)-P(z+2c)=(P(z)+Q(z+c))-(Q(z+c)+P(z+2c))\equiv\chi_1-\chi_2$ and $Q(z)-Q(z+2c)\equiv\chi_2-\chi_1$. It is easy to see that $P(z)=\sum_{j=1}^nb_jz_j+\Psi_1(z)+A$ and $Q(z)=\sum_{j=1}^nd_jz_j+\Psi_2(z)+B$
where $b_j,d_j,A,B\in\mathbb{C}$ $(1\leq j\leq n)$ and $\Psi_k(z)$ $(k=1,2)$ is a polynomial defined in (\ref{K1}). From (\ref{bv1}), we have
\beas
-a_4\frac{\pa^2 \Psi_1(z)}{\pa z_\mu^2}-a_1i\left(b_\mu+\frac{\pa \Psi_1(z)}{\pa z_\mu}\right)-a_4\left(b_\mu+\frac{\pa \Psi_1(z)}{\pa z_\mu}\right)^2\equiv a_3K_4e^{-\chi_1}/K_1+a_2\;\;\text{and}\\
-a_4\frac{\pa^2 \Psi_2(z)}{\pa z_\mu^2}-a_1i\left(b_\mu+\frac{\pa \Psi_2(z)}{\pa z_\mu}\right)-a_4\left(b_\mu+\frac{\pa \Psi_2(z)}{\pa z_\mu}\right)^2\equiv a_3K_2e^{-\chi_2}/K_3+a_2. \eeas
If $\Psi_k(z)$ $(k=1,2)$ contain the variable $z_\mu$, then by comparing the degrees on both sides, we get 
that $\deg(\Psi_k(t))\leq 1$. For simplicity, we still denote $P(z)=\sum_{j=1}^nb_jz_j+A$ and $Q(z)=\sum_{j=1}^nd_jz_j+B$, where $b_j,d_j,A,B\in\mathbb{C}$ $(1\leq j\leq n)$. This implies that $\Psi_k(z)\equiv 0$ for $k=1,2$. 
Since $P(z)+Q(z+c)$ is a constant, thus we have $b_j+d_j=0$ for $1\leq j\leq n$. 
Therefore $P(z)=\sum_{j=1}^nb_jz_j+A$ and $Q(z)=-\sum_{j=1}^nb_jz_j+B$.
From (\ref{bv1}) and (\ref{bv5}), we deduce that
\be\label{bv6}\left\{\begin{array}{ll}
\frac{a_1K_1}{a_3K_4}\left(\frac{1}{i}b_\mu-\frac{a_4}{a_1}b_\mu^2-\frac{a_2}{a_1}\right)e^{-\sum_{j=1}^nb_jc_j+A+B}\equiv 1,\\[1mm]
\frac{a_1K_3}{a_3K_2}\left(-\frac{1}{i}b_\mu-\frac{a_4}{a_1}b_\mu^2-\frac{a_2}{a_1}\right)e^{\sum_{j=1}^nb_jc_j+A+B}\equiv 1,\\[1mm]
\frac{a_1K_2}{a_3K_3}\left(\frac{1}{i}b_\mu-\frac{a_4}{a_1}b_\mu^2-\frac{a_2}{a_1}\right)e^{\sum_{j=1}^nb_jc_j-A-B}\equiv 1,\\[1mm]
\frac{a_1K_4}{a_3K_1}\left(-\frac{1}{i}b_\mu-\frac{a_4}{a_1}b_\mu^2-\frac{a_2}{a_1}\right)e^{-\sum_{j=1}^nb_jc_j-A-B}\equiv 1.\end{array}\right.\ee
From (\ref{bv6}), we have 
\beas &&\left(\frac{a_1}{a_3}\right)^2\left(\frac{1}{i}b_\mu-\frac{a_4}{a_1}b_\mu^2-\frac{a_2}{a_1}\right)^2=\left(\frac{a_1}{a_3}\right)^2\left(-\frac{1}{i}b_\mu-\frac{a_4}{a_1}b_\mu^2-\frac{a_2}{a_1}\right)^2,\\[1mm]\text{\it i.e.,}
&&\frac{1}{i}b_\mu\left(\frac{a_4}{a_1}b_\mu^2+\frac{a_2}{a_1}\right)=0,\;\text{\it i.e.,}\; \text{either}\quad  b_\mu=0\quad\text{or}\quad a_4b_\mu^2+a_2=0.\eeas
It is clear that both $b_\mu$ and $a_4b_\mu^2+a_2$ are not simultaneously zero, otherwise we get $a_2=0$, which is a contradiction. Now two different cases possible which are 
$b_\mu=0$ and $a_4b_\mu^2+a_2=0$. The second case is consider in the proof of \textrm{Theorem \ref{T12}}.\\[2mm]
If $b_\mu=0$, then we deduce from (\ref{bv6}) that 
\be\label{1bx9}\left\{\begin{array}{lll}
e^{2\sum_{j=1,j\not=\mu}^nb_jc_j}\equiv \left(a_3/a_2\right)^2\equiv e^{-2\sum_{j=1,j\not=\mu}^nb_jc_j},\quad e^{2A+2B}\equiv a_3^2K_2K_4/(a_2^2K_1K_3),\\[2mm]
e^{-\sum_{j=1,j\not=\mu}^nb_jc_j+A+B}\equiv -a_3K_4/(a_2K_1), \quad e^{\sum_{j=1,j\not=\mu}^nb_jc_j-A-B}\equiv -a_3K_3/(a_2K_2),\\[2mm]
e^{\sum_{j=1,j\not=\mu}^nb_jc_j+A+B}\equiv -a_3K_2/(a_2K_3),  \quad e^{-\sum_{j=1,j\not=\mu}^nb_jc_j-A-B}\equiv -a_3K_1/(a_2K_4).\end{array}\right.\ee
From (\ref{1bx9}), it is clear that $a_2^2=\pm a_3^2$.
From (\ref{eq1}), we deduce that
 \bea\label{bx1} \begin{array}{lll}
 f_1(z)=\frac{z_\mu}{2a_1}\left(K_1e^{\sum_{j=1,j\not=\mu}^nb_jz_j+A}+K_2e^{-\sum_{j=1,j\not=\mu}^nb_jz_j-A}\right)+h_3(y_1),\\[2mm]
 f_2(z)=\frac{z_\mu}{2a_1}\left(K_3e^{-\sum_{j=1,j\not=\mu}^nb_jz_j+B}+K_4e^{\sum_{j=1,j\not=\mu}^nb_jz_j-B}\right)+h_4(y_1),\end{array}\eea
where $h_j(y_1)$ $(j=3,4)$ are finite order entire functions. 
By utilizing the equations  (\ref{1bx9}) and (\ref{bx1}) in the second and fourth equations of (\ref{eq1}), we deduce that respectively
\beas &&a_2h_3(y_1)+a_3h_4(y_1+s_1)=\gamma_1(1)K_1e^{\sum_{j=1,j\not=\mu}^nb_jz_j+A}+\gamma_2(1)K_2e^{-\sum_{j=1,j\not=\mu}^nb_jz_j-A}\\[1mm]\text{and}
&&a_2h_4(y_1)+a_3h_3(y_1+s_1)=\gamma_1(1)K_3e^{-\sum_{j=1,j\not=\mu}^nb_jz_j+B}+\gamma_2(1)K_4e^{\sum_{j=1,j\not=\mu}^nb_jz_j-B},\eeas
where $\gamma_1(1)$ and $\gamma_2(1)$ are given in (\ref{K2}).\\[1.5mm]
If $\Psi_k(z)$ $(k=1,2)$ does not contain the variables $z_\mu$, then, we have $P(z)=\sum_{j=1}^nb_jz_j+\Psi_1(z)+A$ and $Q(z)=\sum_{j=1}^nd_jz_j+\Psi_2(z)+B$, where $b_j,d_j,A,B\in\mathbb{C}$ $(1\leq j\leq n)$ and $\Psi_k(z)$ $(k=1,2)$ is a polynomial defined in (\ref{K1}).
Since $P(z)+Q(z+c)$ is a constant, so we must have $b_j+d_j=0$ for $1\leq j\leq n$ and $\Psi_1(z)+\Psi_2(z)\equiv 0$.
Therefore $P(z)=\sum_{j=1}^nb_jz_j+\Psi_1(z)+A$ and $Q(z)=-\sum_{j=1}^nb_jz_j-\Psi_1(z)+B$.
From (\ref{bv1}) and (\ref{bv5}), we again have (\ref{bv6}) and either $b_\mu=0$ or $a_4b_\mu^2+a_2=0$.  The second case is consider in the proof of \textrm{Theorem \ref{T12}}\\[1.5mm]
If $b_\mu=0$, then we have (\ref{1bx9}). By means of argument similar to the ones above, we have  
\beas \begin{array}{lll} f_1(z)=\frac{z_\mu}{2a_1}\left(K_1e^{\sum_{j=1,j\not=\mu}^nb_jz_j+\Psi_1(z)+A}+K_2e^{-\sum_{j=1,j\not=\mu}^nb_jz_j-\Psi_1(z)-A}\right)+g_1(y_1),\\[2mm]
 f_2(z)=\frac{z_\mu}{2a_1}\left(K_3e^{-\sum_{j=1,j\not=\mu}^nb_jz_j-\Psi_1(z)+B}+K_4e^{\sum_{j=1,j\not=\mu}^nb_jz_j+\Psi_1(z)-B}\right)+g_2(y_1),\end{array}\eeas
where $a_2^2=\pm a_3^2$ and $g_j(y_1)$ $(j=1,2)$ are finite order entire functions satisfying 
\beas a_2g_1(y_1)+a_3g_2(y_1+s_1)=\gamma_1(1)K_1e^{\sum_{j=1,j\not=\mu}^nb_jz_j+A}+\gamma_2(1)K_2e^{-\sum_{j=1,j\not=\mu}^nb_jz_j-A}\\[1mm]\text{and}\;
a_2g_2(y_1)+a_3g_1(y_1+s_1)=\gamma_1(1)K_3e^{-\sum_{j=1,j\not=\mu}^nb_jz_j+B}+\gamma_2(1)K_4e^{\sum_{j=1,j\not=\mu}^nb_jz_j-B},\eeas
where $\gamma_1(1)$ and $\gamma_2(1)$ are given in (\ref{K2}).\end{proof}
\begin{proof}[\bf{Proof of Theorem \ref{T12}}] Let $\Psi_k(z)$ $(k=1,2)$ contain the variable $z_\mu$ and considering the proof of \textrm{Theorem \ref{T11}}, we assume that $a_4b_\mu^2+a_2=0$. From (\ref{bv6}) we have
\be\label{bv9}\left\{\begin{array}{lll}
e^{2\sum_{j=1}^nb_jc_j}\equiv \left(a_3/(a_1b_\mu)\right)^2\equiv e^{-2\sum_{j=1}^nb_jc_j}, \quad e^{2A+2B}\equiv a_3^2K_2K_4/(a_1^2b_\mu^2K_1K_3)\\[1.5mm]
e^{-\sum_{j=1}^nb_jc_j+A+B}\equiv ia_3K_4/(a_1b_\mu K_1), \quad e^{\sum_{j=1}^nb_jc_j-A-B}\equiv ia_3K_3/(a_1b_\mu K_2),\\[1.5mm]
e^{\sum_{j=1}^nb_jc_j+A+B}\equiv -ia_3K_2/(a_1b_\mu K_3), \quad e^{-\sum_{j=1}^nb_jc_j-A-B}\equiv -ia_3K_1/(a_1b_\mu K_4).\end{array}\right.\ee
From (\ref{bv9}), we have $a_3^2=\pm (a_1b_\mu)^2$.
The Lagrange's auxiliary equations \cite[Chapter 2]{430} of the first equation of (\ref{eq1}) are 
\beas \frac{dz_1}{0}=\cdots=\frac{dz_{\mu-1}}{0}=\frac{dz_\mu}{1}=\frac{dz_{\mu+1}}{0}=\cdots=\frac{dz_n}{0}=\frac{2a_1df_1(z)}{K_1e^{\sum_{j=1}^nb_jz_j+A}+K_2e^{-\sum_{j=1}^nb_jz_j-A}}.\eeas
Note that $\alpha_j=z_j$ for $1\leq j(\not=\mu)\leq n$ and 
\beas&& df_1(z)=\frac{K_1}{2a_1}e^{\sum_{j=1}^nb_jz_j+A}dz_\mu+\frac{K_2}{2a_1}e^{-\sum_{j=1}^nb_jz_j-A}dz_\mu,\\\text{\it{i.e.,}}
&&df_1(z)=\frac{K_1}{2a_1}e^{b_\mu z_\mu+\sum_{j=1,j\not=\mu}^nb_j\alpha_j+A}dz_1+\frac{K_2}{2a_1}e^{-b_\mu z_\mu-\sum_{j=1,j\not=\mu}^nb_j\alpha_j-A}dz_\mu,\eeas
which implies
\beas f_1(z)=\frac{K_1}{2a_1b_\mu}e^{\sum_{j=1}^nb_jz_j+A}-\frac{K_2}{2a_1b_\mu}e^{-\sum_{j=1}^nb_jz_j-A}+\alpha_\mu.\eeas 
Note that after integration with respect to 
$z_\mu$, replacing $\alpha_j$ by $z_j$ $(1\leq j(\not=\mu)\leq n)$, where $\alpha_j\in\mathbb{C}$ for $1\leq j\leq n$. Hence the solution is $\Phi(\alpha_1,\alpha_2,\ldots, \alpha_n )=0$. For simplicity, we suppose
\bea\label{bv7} f_1(z)=\frac{K_1}{2a_1b_\mu}e^{\sum_{j=1}^nb_jz_j+A}-\frac{K_2}{2a_1b\mu}e^{-\sum_{j=1}^nb_jz_j-A}+h_5(y_1),\eea
where $h_5(y_1)$ is a finite order entire function.
Similarly, from the third equation of (\ref{eq1}), we have
\bea\label{bv8} f_2(z)=-\frac{K_3}{2a_1b_\mu}e^{-\sum_{j=1}^nb_jz_j+B}+\frac{K_4}{2a_1b_\mu}e^{\sum_{j=1}^nb_jz_j-B}+h_6(y_1),\eea
where $h_6(y_1)$ is a finite order entire function. 
By utilizing the equations (\ref{bv9}), (\ref{bv7}) and (\ref{bv8}) in the second and fourth equations of (\ref{eq1}), we deduce respectively that
$a_2h_5(y_1)+a_3h_6(y_1+s_1)\equiv0$ and $a_2h_6(y_1)+a_3h_5(y_1+s_1)\equiv0$.\\[2mm]
Let $\Psi_k(z)$ $(k=1,2)$ does not contain the variables $z_\mu$ and considering the proof of \textrm{Theorem \ref{T11}}, we assume that $a_4b_\mu^2+a_2=0$.
Then, again we have (\ref{bv9}). Using argument similar to those presented above, we deduce from (\ref{eq1}) that
\beas \begin{array}{lll}
f_1(z)=\frac{K_1}{2a_1b_\mu}e^{\sum_{j=1}^nb_jz_j+\Psi_1(z)+A}-\frac{K_2}{2a_1b\mu}e^{-\sum_{j=1}^nb_jz_j-\Psi_1(z)-A}+g_3(y_1),\\[2mm]
f_2(z)=-\frac{K_3}{2a_1b_\mu}e^{-\sum_{j=1}^nb_jz_j-\Psi_1(z)+B}+\frac{K_4}{2a_1b_\mu}e^{\sum_{j=1}^nb_jz_j+\Psi_1(z)-B}+g_4(y_1),\end{array}\eeas
where $a_3^2=\pm (a_1b_\mu)^2$ and $g_j(y_1)$ $(j=3,4)$ are finite order entire functions satisfying 
\beas a_2g_3(y_1)+a_3g_4(y_1+s_1)\equiv0\quad\text{and}\quad a_2g_4(y_1)+a_3g_3(y_1+s_1)\equiv0.\eeas \end{proof}
\begin{proof}[\bf{Proof of Theorem \ref{T13}}] 
 Let \bea\label{bv2}\begin{array}{lll}
\Gamma_1(z)e^{P(z)+Q(z+c)}\equiv 1\;\;\text{and}\;\;\Psi_2(z) e^{-Q(z)+P(z+c)}\equiv 1.\end{array}\eea
Since $P(z),Q(z)$ are non-constant polynomials, from (\ref{bv2}) it is evident that $P(z)+Q(z+c)\equiv \xi_1$ and $-Q(z)+P(z+c)\equiv \xi_2$, where $\xi_1,\xi_2\in\mathbb{C}$. Therefore, we have
$P(z)+P(z+2c)\equiv \xi_1+\xi_2$, which is a contradiction.\\
Let \bea\label{bv3}\begin{array}{lll}
\Psi_1(z) e^{-P(z)+Q(z+c)}\equiv1\;\;\text{and}\;\;\Gamma_2(z)e^{Q(z)+P(z+c)}\equiv 1.\end{array}\eea
From (\ref{bv3}), it is clear that $-P(z)+Q(z+c)\equiv \xi_3$ and $Q(z)+P(z+c)\equiv \xi_4$, where $\xi_3,\xi_4\in\mathbb{C}$. 
Then, by using similar argument as above, we get a contradiction.\\
 Let \bea\label{bv4}\begin{array}{lll}
\Psi_1(z) e^{-P(z)+Q(z+c)}\equiv1\;\;\text{and}\;\;\Psi_2(z) e^{-Q(z)+P(z+c)}\equiv 1.\end{array}\eea
Using (\ref{bv4}), we deduce from (\ref{tr}) and (\ref{14p10}) that
\bea\label{1bv5}\begin{array}{lll}
K_4\Gamma_1(z) e^{P(z)-Q(z+c)}/K_3\equiv1\quad\text{and}\quad K_2\Gamma_2(z) e^{Q(z)-P(z+c)}/K_1\equiv 1.\end{array}\eea
From (\ref{bv4}), it is clear that $P(z)-Q(z+c)$ and $Q(z)-P(z+c)$ are both constants, say $\chi_1$ and $\chi_2$ respectively, where $\chi_1,\chi_2\in\mathbb{C}$. Now, 
$P(z)-P(z+2c)=(P(z)-Q(z+c))+(Q(z+c)-P(z+2c))\equiv\chi_1+\chi_2$ and $Q(z)-Q(z+2c)\equiv\chi_1+\chi_2$. Therefore,
$P(z)=\sum_{j=1}^nb_jz_j+\Psi_1(z)+A$ and $Q(z)=\sum_{j=1}^nd_jz_j+\Psi_2(z)+B$, where $b_j,d_j,A,B\in\mathbb{C}$ for $1\leq j\leq n$ and and $\Psi_k(z)$ $(k=1,2)$ is a polynomial defined in (\ref{K1}). From (\ref{bv4}), we have 
\bea\label{K3}\begin{array}{lll}
\qquad-i\left(b_\mu+\frac{\pa\Psi_1(z)}{\pa z_\mu}\right)-\frac{a_4}{a_1}\left(b_\mu+\frac{\pa\Psi_1(z)}{\pa z_\mu}\right)^2+\frac{a_4}{a_1}\frac{\pa^2\Psi_1(z)}{\pa z_\mu^2}-\frac{a_2}{a_1}\equiv \frac{a_3K_4}{a_1K_2}e^{\chi_1}\\[2mm]\text{and}\quad
-i\left(d_\mu+\frac{\pa\Psi_2(z)}{\pa z_\mu}\right)-\frac{a_4}{a_1}\left(d_\mu+\frac{\pa\Psi_2(z)}{\pa z_\mu}\right)^2+\frac{a_4}{a_1}\frac{\pa^2\Psi_2(z)}{\pa z_\mu^2}-\frac{a_2}{a_1}\equiv \frac{a_3K_2}{a_1K_4}e^{\chi_2}.\end{array}\eea
{\bf Case 1.} If $\Psi_k(z)$ $(k=1,2)$ contain the variable $z_\mu$, then by comparing the degrees on both sides of (\ref{K3}), we get that $\deg(\Psi_k(z))\leq1$ for $k=1,2$. For simplicity, we still denote $P(z)=\sum_{j=1}^nb_jz_j+A$ and $Q(z)=\sum_{j=1}^nd_jz_j+B$, where $b_j,d_j\in\mathbb{C}$ $(1\leq j\leq n)$.  This implies that $\Psi_k (z)\equiv 0$ for $k=1,2$. 
Since $P(z)-Q(z+c)$ is a constant, so we must have $b_j=d_j$ for $1\leq j\leq n$.
Therefore $P(z)=\sum_{j=1}^nb_jz_j+A$ and $Q(z)=\sum_{j=1}^nb_jz_j+B$, where $b_j,A,B\in\mathbb{C}$ for $1\leq  j\leq n$.
From (\ref{bv4}) and (\ref{1bv5}), we obtain
\be\label{1bv6}\left\{\begin{array}{ll}
\frac{a_1K_2}{a_3K_4}\left(\frac{1}{i}b_\mu-\frac{a_4}{a_1}b_\mu^2-\frac{a_2}{a_1}\right)e^{\sum_{j=1}^nb_jc_j-A+B}\equiv 1,\\[1mm]
\frac{a_1K_4}{a_3K_2}\left(\frac{1}{i}b_\mu-\frac{a_4}{a_1}b_\mu^2-\frac{a_2}{a_1}\right)e^{\sum_{j=1}^nb_jc_j+A-B}\equiv 1,\\[1mm]
\frac{a_1K_1}{a_3K_3}\left(\frac{1}{i}b_\mu-\frac{a_4}{a_1}b_\mu^2-\frac{a_2}{a_1}\right)e^{-\sum_{j=1}^nb_jc_j+A-B}\equiv 1,\\[1mm]
\frac{a_1K_3}{a_3K_1}\left(\frac{1}{i}b_\mu-\frac{a_4}{a_1}b_\mu^2-\frac{a_2}{a_1}\right)e^{-\sum_{j=1}^nb_jc_j-A+B}\equiv 1.\end{array}\right.\ee
From (\ref{1bv6}), we have 
\beas \left(-i b_\mu-a_4b_\mu^2/a_1-a_2/a_1\right)^2=\left(a_3/a_1\right)^2,
\quad\text{\it{i.e.,}}\quad
ia_1b_\mu+a_4b_\mu^2+a_2=\pm a_3.\eeas
If $ia_1b_\mu+a_4b_\mu^2+a_2=(-1)^\nu a_3$ $(\nu=1,2)$, then by similar argument to those of \textrm{Theorem \ref{T12}}, we obtain 
\bea\label{K6}\begin{array}{lll} 
f_1(z)=\frac{K_1}{2a_1b_\mu}e^{\sum_{j=1}^nb_jz_j+A}-\frac{K_2}{2a_1b_\mu}e^{-\sum_{j=1}^nb_jz_j-A}+h_7(y_1),\\[2mm]
f_2(z)=\frac{K_3}{2a_1b_\mu}e^{\sum_{j=1}^nb_jz_j+B}-\frac{K_4}{2a_1b_\mu}e^{-\sum_{j=1}^nb_jz_j-B}+h_8(y_1),\end{array}\eea
where $b_j,A,B\in\mathbb{C}$ $(1\leq j\leq n )$, $h_j(y_1)$ $( j=7, 8)$ are finite order entire functions satisfying $a_2h_7(y_1)+a_3h_8(y_1+s_1)\equiv0$ and $a_2h_8(y_1)+a_3h_7(y_1+s_1)\equiv0$ with
\bea\label{K4}\left\{\begin{array}{llll}
e^{2\sum_{j=1}^nb_jc_j}\equiv 1, e^{2A-2B}\equiv K_2K_3/(K_1K_4), e^{\sum_{j=1}^nb_jc_j-A+B}\equiv (-1)^{\nu+1}K_4/K_2,\\[1.5mm]
 e^{\sum_{j=1}^nb_jc_j+A-B}\equiv (-1)^{\nu+1}K_2/K_4, e^{-\sum_{j=1}^nb_jc_j+A-B}\equiv (-1)^{\nu+1}K_3/K_1,\\[1.5mm] 
 e^{-\sum_{j=1}^nb_jc_j-A+B}\equiv (-1)^{\nu+1}K_1/K_3.\end{array}\right.\eea
{\bf Case 2.} If $\Psi_k(z)$ $(k=1,2)$ does not contain the variable $z_\mu$, then, we have $P(z)=\sum_{j=1}^nb_jz_j+\Psi_1(z)+A$ and $Q(z)=\sum_{j=1}^nd_jz_j+\Psi_2(z)+B$, where $b_j,d_j,A,B\in\mathbb{C}$ $(1\leq j\leq n)$ and $\Psi_k(z)$ $(k=1,2)$ is a polynomial defined in (\ref{K1}).
Since $P(z)-Q(z+c)$ is a constant, so we must have $b_j=d_j$ for $1\leq j\leq n$ and $\Psi_1(z)\equiv \Psi_2(z)$.
Therefore $P(z)=\sum_{j=1}^nb_jz_j+\Psi_1(z)+A$ and $Q(z)=\sum_{j=1}^nb_jz_j+\Psi_1(z)+B$, where $b_j,A,B\in\mathbb{C}$ for $1\leq j\leq n$. From (\ref{bv4}) and (\ref{1bv5}), we get (\ref{1bv6}) and $ia_1b_\mu+a_4b_\mu^2+a_2=\pm a_3$.
If $ia_1b_\mu+a_4b_\mu^2+a_2=(-1)^\nu a_3$ $(\nu=1,2)$, then by similar argument to those of \textrm{Case 1}, we obtain 
\bea\label{K7}\begin{array}{lll}
f_1(z)=\frac{K_1}{2a_1b_\mu}e^{\sum_{j=1}^nb_jz_j+\Psi_1(z)+A}-\frac{K_2}{2a_1b_\mu}e^{-\sum_{j=1}^nb_jz_j-\Psi_1(z)-A}+g_5(y_1),\\[3mm]
f_2(z)=\frac{K_3}{2a_1b_\mu}e^{\sum_{j=1}^nb_jz_j+\Psi_1(z)+B}-\frac{K_4}{2a_1b_\mu}e^{-\sum_{j=1}^nb_jz_j-\Psi_1(z)-B}+g_6(y_1),\end{array}\eea
where $b_j,A,B\in\mathbb{C}$ $(1\leq j\leq n )$, $g_j(y_1)$ $( j=5,6)$ are finite order entire functions satisfying $a_2g_5(y_1)+a_3g_6(y_1+s_1)\equiv0$ and $a_2g_6(y_1)+a_3g_5(y_1+s_1)\equiv0$ with (\ref{K4}).
This completes the proof.
\end{proof}
\begin{proof}[\bf{Proof of Theorem \ref{T2}}] {\bf Part 1.} The first part is common for all Theorems \ref{T2}-\ref{T24}.
Let $(f_1,f_2)$ be a pair of finite order transcendental entire functions  satisfies the system (\ref{e4}). By using similar argument to those of \textrm{Theorem \ref{T1}}, we get
\bea\label{teq1}\left\{\begin{array}{llll}
a_1\frac{\pa f_1(z)}{\pa z_1}+a_2\frac{\pa f_1(z)}{\pa z_2}+\cdots+a_n\frac{\pa f_1(z)}{\pa z_n}=\frac{K_1e^{P(z)}+K_2e^{-P(z)}}{2},\\[2mm]
a_{n+1}f_1(z)+a_{n+2}f_2(z+c)=\frac{K_1e^{P(z)}-K_2e^{-P(z)}}{2i},\\[2mm]
a_1\frac{\pa f_2(z)}{\pa z_1}+a_2\frac{\pa f_2(z)}{\pa z_2}+\cdots+a_n\frac{\pa f_2(z)}{\pa z_n}=\frac{K_3e^{Q(z)}+K_4e^{-Q(z)}}{2},\\[2mm]
a_{n+1}f_2(z)+a_{n+2}f_1(z+c)=\frac{K_3e^{Q(z)}-K_4e^{-Q(z)}}{2i},\end{array}\right.\eea
where $K_1,K_2,K_3,K_4\in\mathbb{C}\setminus\{0\}$ such that $K_1K_2=1=K_3K_4$ and $P(z)$ and $Q(z)$ are polynomials on $\mathbb{C}^n$. The following cases arise separately in the proofs of all \textrm{Theorems} \ref{T2}-\ref{T24}.\\
{\bf Part 2.} We now begin to prove \textrm{Theorem \ref{T2}} properly.  Let $P(z)$ and $Q(z)$ be both constants. From (\ref{teq1}), we have 
\bea\label{J1}
&&\sum_{j=1}^n a_j\frac{\pa f_1(z)}{\pa z_j}=A_1,\quad a_{n+1}f_1(z)+a_{n+2}f_2(z+c)=A_2,\\[1.5mm]
\label{J2}&&\sum_{j=1}^n a_j\frac{\pa f_2(z)}{\pa z_j}=A_3\quad\text{and }\quad a_{n+1}f_2(z)+a_{n+2}f_1(z+c)=A_4,\eea
where $A_j\in\mathbb{C}$ for $1\leq j\leq 4$ with $A_k^2+A_{k+1}^2=1$ $(k=1,3)$.
The Lagrange's auxiliary equations \cite[Chapter 2]{430} of the first equation in (\ref{J1}) are
\beas \frac{dz_1}{a_1}=\frac{dz_2}{a_2}=\frac{dz_3}{a_3}=\cdots=\frac{dz_n}{a_n}=\frac{df_1(z)}{A_1}.\eeas
Note that $z_j=(\beta_j+a_jz_\mu)/a_\mu$ for $1\leq j(\not=\mu)\leq n$, where $\mu\in\{1,2,\ldots,n\}$ and $df_1(z)=(A_1/a_\mu)dz_\mu$, {\it i.e.,} $f_1(z)=(A_1/a_\mu)z_\mu+\beta_1$, where $\beta_j\in\mathbb{C}$ for $1\leq j\leq n$.
Hence the solution is $\Phi(\beta_1,\beta_2,\ldots, \beta_n )=0$. For simplicity, we suppose
\bea\label{req2} f_1(z)=(A_1/a_\mu)z_\mu+g_1\left(y\right),\eea
where $g_1\left(y\right)$ is a finite order transcendental entire function of $a_\mu z_1-a_1z_\mu, a_\mu z_2-a_2z_\mu,\ldots, a_\mu z_{\mu-1}-a_{\mu-1}z_\mu, a_\mu z_{\mu+1}-a_{\mu+1}z_\mu,\ldots,a_\mu z_n-a_n z_\mu$ such that $\sum_{j=1}^n a_j\frac{\pa g_1(z)}{\pa z_j}\equiv 0$.
From the first equation of (\ref{J2}), we obtain by using similar argument as above
\bea\label{req3}  f_2(z)=(A_3/a_\mu)z_\mu+g_2(y),\eea
where $g_2(y)$ is a finite order transcendental entire function with $\sum_{j=1}^n a_j\frac{\pa g_2(z)}{\pa z_j}\equiv 0$. Using (\ref{req2}) and (\ref{req3}), we get from the second equation of (\ref{J1}) that
\be\label{mnb}(a_{n+1} A_1+a_{n+2}A_3)z_\mu/a_\mu +a_{n+1}g_1(y)+a_{n+2}g_2(y+s)+(a_{n+2}c_\mu A_3)/a_\mu\equiv A_2.\ee
Comparing both sides of (\ref{mnb}), we get
\beas a_{n+1}A_1+a_{n+2}A_3=0,\;a_{n+1}g_1(y)+a_{n+2}g_2(y+s)\equiv 0\quad\text{and}\quad a_{n+2}c_\mu A_3= a_\mu A_2.\eeas
Similarly, by using (\ref{req2}) and (\ref{req3}), we get from the second equation of (\ref{J2}) that
\beas a_{n+1}A_3+a_{n+2}A_1= 0,\;a_{n+1}g_2(y)+a_{n+2}g_1(y+s)\equiv 0\quad\text{and}\quad a_{n+2}c_\mu A_1= a_\mu A_4.\eeas
Using similar argument as in \textrm{Case 1} of \textrm{Theorem \ref{T1}}, we have $a_{n+1}=\pm a_{n+2}$, 
\beas f_1(z)=\frac{z_\mu}{\sqrt{a_\mu^2+a_{n+1}^2c_\mu^2}}+g_1(y)\;\;\text{and}\;\; f_2(z)=\frac{z_\mu}{\sqrt{a_\mu^2+a_{n+1}^2c_\mu^2}}+g_2(y),\eeas
where $1\leq \mu\leq n$ and $g_j(y)$ $(j=1,2)$ are finite order transcendental entire functions with periods $2s$ satisfying $\sum_{k=1}^n a_k\frac{\pa g_j(y)}{\pa z_k}\equiv 0$.\end{proof}
\begin{proof}[\bf{Proof of Theorem \ref{T22}}]
Let either $P(z)$ or $Q(z)$ be a constant. By using similar argument as of \textrm{Theorem \ref{T11}}, we get a contradiction.\\
 Let $P(z),Q(z)$ be both non-constant polynomials. Differentiating partially with respect to $z_j$ on both sides of the second equation in (\ref{teq1}), we deduce that
\bs\bea\label{mnb2}a_{n+1}\frac{\pa f_1(z)}{\pa z_j}+a_{n+2}\frac{\pa f_2(z+c)}{\pa z_j}=\frac{K_1e^{P(z)}+K_2e^{-P(z)}}{2i}\frac{\pa P(z)}{\pa z_j}\;(1\leq j\leq n).\eea\es
From (\ref{teq1}) and (\ref{mnb2}), we can conclude that
\bea
\label{mnb4}\frac{K_1}{ia_{n+2}K_4}\Gamma_3(z)e^{P(z)+Q(z+c)}-\frac{K_3}{K_4}e^{2Q(z+c)}+\frac{K_2}{ia_{n+2}K_4}\Gamma_3(z) e^{-P(z)+Q(z+c)}\equiv 1,\eea
where $\Gamma_3(z)=\sum_{j=1}^n a_j\frac{\pa P(z)}{\pa z_j}-ia_{n+1}$. Similarly, from (\ref{teq1}), we get
\bea\label{mnb5}\frac{K_3}{ia_{n+2}K_2}\Gamma_4(z)e^{Q(z)+P(z+c)}-\frac{K_1}{K_2}e^{2P(z+c)}+\frac{K_4}{ia_{n+2}K_2}\Gamma_4(z)e^{-Q(z)+P(z+c)}\equiv 1,\eea
where $\Gamma_4(z)=\sum_{j=1}^n a_j\frac{\pa Q(z)}{\pa z_j}-ia_{n+1}$.
From (\ref{mnb4}) and (\ref{mnb5}), we observe that $\Gamma_3(z)\not\equiv 0$ and $\Gamma_4(z)\not\equiv 0$, otherwise we get a contradiction. 
From (\ref{mnb4}), it is evident that
\bs\beas &&N\left(r, \frac{K_1}{ia_{n+2}K_4}\Gamma_3(z)e^{P(z)+Q(z+c)}\right)=N\left(r, 0;\frac{K_1}{ia_{n+2}K_4}\Gamma_3(z)e^{P(z)+Q(z+c)}\right)\\
&&=N\left(r, \frac{K_2}{ia_{n+2}K_4}\Gamma_3(z)e^{-P(z)+Q(z+c)}\right)=N\left(r, 0; \frac{K_2}{ia_{n+2}K_4}\Gamma_3(z)e^{-P(z)+Q(z+c)}\right)\\
&&=N\left(r, -\frac{K_3}{K_4}e^{2Q(z+c)}\right)=N\left(r, 0;-\frac{K_3}{K_4}e^{2Q(z+c)}\right)=S\left(r,-\frac{K_3}{K_4}e^{2Q(z+c)}\right).\eeas\es
In view of \textrm{Lemma \ref{lem1}} and from (\ref{mnb4}), we have 
\beas\text{either}\quad \frac{K_1}{ia_{n+2}K_4}\Gamma_3(z)e^{P(z)+Q(z+c)}\equiv 1\quad \text{or}\quad \frac{K_2}{ia_{n+2}K_4}\Gamma_3(z)e^{-P(z)+Q(z+c)}\equiv1.\eeas
 In view of \textrm{Lemma \ref{lem1}} and by using similar argument as above, we get from (\ref{mnb5}) that 
\beas\text{either}\quad \frac{K_3}{ia_{n+2}K_2}\Gamma_4(z)e^{Q(z)+P(z+c)}\equiv 1\quad\text{or}\quad \frac{K_4}{ia_{n+2}K_2}\Gamma_4(z)e^{-Q(z)+P(z+c)}\equiv 1.\eeas
Now we will discuss the following cases.\\ 
Let 
\bea\label{w1}\begin{array}{ll}
\frac{K_1}{ia_{n+2}K_4}\Gamma_3(z)e^{P(z)+Q(z+c)}\equiv 1\quad \text{and}\quad
\frac{K_3}{ia_{n+2}K_2}\Gamma_4(z)e^{Q(z)+P(z+c)}\equiv 1.\end{array}\eea
Using (\ref{w1}), we deduce from (\ref{mnb4}) and (\ref{mnb5}) respectively
 \bea\label{w11}\begin{array}{ll}
 \frac{K_2}{ia_{n+2}K_3}\Gamma_3(z)e^{-P(z)-Q(z+c)}\equiv 1\quad \text{and}\quad
\frac{K_4}{ia_{n+2}K_1}\Gamma_4(z)e^{-Q(z)-P(z+c)}\equiv 1.\end{array}\eea
It is clear that $P(z)+Q(z+c)\equiv \chi_1$ and $Q(z)+P(z+c)\equiv \chi_2$, where $\chi_1,\chi_2\in\mathbb{C}$.
By using similar argument to those of \textrm{Sub-case 3.1} of \textrm{Theorem \ref{T1}}, we have
$P(z)=\sum_{j=1}^nb_jz_j+\Psi_1(z)+A$ and $Q(z)=\sum_{j=1}^nd_jz_j+\Psi_2(z)+d_{n+1}$, where $b_j,d_j,A,B\in\mathbb{C}$ $(1\leq j\leq n)$ and $\Psi_k(z)$ $(k=1,2)$ is a polynomial defined in (\ref{K1}). From (\ref{w1}), we have
\be\label{K8}\left\{\begin{array}{ll}\sum_{j=1}^n a_j\left(b_j+\frac{\pa\Phi_1(z)}{\pa z_j}\right)-ia_{n+1}\equiv \frac{ia_{n+2}K_4}{K_1}e^{-\chi_1},\\
\sum_{j=1}^n a_j\left(d_j+\frac{\pa\Phi_2(z)}{\pa z_j}\right)-ia_{n+1}\equiv \frac{ia_{n+2}K_2}{K_3}e^{-\chi_2}.\end{array}\right.\ee
Since $a_j\not=0$ for $j=1,2,\ldots,n$, by comparing the degrees on both sides of (\ref{K8}), we get that $\deg(\Psi_k(z))\leq1$ for $k=1,2$. For simplicity, we still denote 
$P(z)=\sum_{j=1}^nb_jz_j+A$ and $Q(z)=\sum_{j=1}^nd_jz_j+B$, where $b_j,d_j,A,B\in\mathbb{C}$ $(1\leq j\leq n)$.  This implies that $\Psi_k (z)\equiv 0$ for $k=1,2$. 
Since $P(z)+Q(z+c)$ is a constant, so we must have $b_j+d_j=0$ for $1\leq j\leq n$.
Therefore $P(z)=\sum_{j=1}^nb_jz_j+A$ and $Q(z)=-\sum_{j=1}^nb_jz_j+B$, where $b_j,A,B\in\mathbb{C}$ for $1\leq j\leq n$.
From (\ref{w1}) and (\ref{w11}), we have
\bea\label{w12}\left\{\begin{array}{ll}
\frac{K_1}{ia_{n+2}K_4}\left(a_1b_1+a_2b_2+\cdots+a_nb_n-ia_{n+1}\right)e^{-\sum_{j=1}^nb_jc_j+A+B}\equiv 1,\\[2mm]
\frac{K_3}{ia_{n+2}K_2}\left(-a_1b_1-a_2b_2\cdots-a_nb_n-ia_{n+1}\right)e^{\sum_{j=1}^nb_jc_j+A+B}\equiv 1,\\[2mm]
 \frac{K_2}{ia_{n+2}K_3}\left(a_1b_1+a_2b_2+\cdots+a_nb_n-ia_{n+1}\right)e^{\sum_{j=1}^nb_jc_j-A-B}\equiv 1,\\[2mm]
\frac{K_4}{ia_{n+2}K_1}\left(-a_1b_1-a_2b_2-\cdots-a_nb_n-ia_{n+1}\right)e^{-\sum_{j=1}^nb_jc_j-A-B}\equiv 1.\end{array}\right.\eea
From (\ref{w12}), we deduce that
\beas &&\left(a_1b_1+a_2b_2+\cdots+a_nb_n-ia_{n+1}\right)^2=\left(-a_1b_1-a_2b_2-\cdots-a_nb_n-ia_{n+1}\right)^2,\\ {\it i.e.}\;
&& ia_{n+1}(a_1b_1+a_2b_2+\cdots+a_nb_n)=0,\; {\it i.e.}, \;a_1b_1+a_2b_2+\cdots+a_nb_n=0.\eeas
From (\ref{w12}), we have
\be\label{1w1}\left\{\begin{array}{lll}
e^{2\sum_{j=1}^nb_jc_j}\equiv \left(a_{n+2}/a_{n+1}\right)^2\equiv e^{-2\sum_{j=1}^nb_jc_j},\quad e^{2A+2B }\equiv a_{n+2}^2K_2K_4/(a_{n+1}^2K_1K_3),\\[1.5mm]
e^{-\sum_{j=1}^nb_jc_j+A+B}\equiv -a_{n+2}K_4/(a_{n+1}K_1), \quad e^{\sum_{j=1}^nb_jc_j+A+B}\equiv -a_{n+2}K_2/(a_{n+1}K_3),\\[1.5mm]
e^{\sum_{j=1}^nb_jc_j-A-B}\equiv -a_{n+2}K_3/(a_{n+1}K_2), \quad e^{-\sum_{j=1}^nb_jc_j-A-B}\equiv -a_{n+2}K_1/(a_{n+1}K_4).\end{array}\right.\ee
From (\ref{w12}), it is clear that $a_{n+1}^2=\pm a_{n+2}^2$.
The Lagrange's auxiliary equations \cite[Chapter 2]{430} of the first equation of (\ref{teq1}) are
\beas \frac{dz_1}{a_1}=\frac{dz_2}{a_2}=\frac{dz_3}{a_3}=\cdots=\frac{dz_n}{a_n}=\frac{2df_1(z)}{K_1e^{\sum_{j=1}^nb_jz_j+A}+K_2e^{-\sum_{j=1}^nb_jz_j-A}}.\eeas
Note that $z_j=(\beta_j+a_jz_\mu)/a_\mu$ for $1\leq j(\not=\mu)\leq n$, where $\mu\in\{1,2,\ldots,n\}$ and 
\beas\sum_{j=1}^nb_jz_j&=&b_\mu z_\mu+\sum_{j=1,j\not=\mu}^n b_j(\beta_j+a_j z_\mu)/a_\mu\\
&=&\left(b_\mu+\sum_{j=1,j\not=\mu}^n b_ja_j/a_\mu\right)z_\mu+\sum_{j=1,j\not=\mu}^n b_j \beta_j/a_\mu\\
&=&(1/a_\mu)\sum_{j=1,j\not=\mu}^n b_j\beta_j.\eeas
Thus,
\beas &&df_1(z)=\frac{K_1}{2a_\mu}e^{\sum_{j=1}^nb_jz_j+A}dz_\mu+\frac{K_2}{2a_\mu}e^{-\sum_{j=1}^nb_jz_j-A}dz_\mu,\\\text{\it{i.e.,}}
&& df_1(z)=\frac{K_1}{2a_\mu}e^{(1/a_\mu)\sum_{j=1,j\not=\mu}^n b_j\beta_j+A}dz_\mu+\frac{K_2}{2a_\mu}e^{-(1/a_\mu)\sum_{j=1,j\not=\mu}^n b_j\beta_j-A}dz_\mu,\eeas
which implies 
\beas f_1(z)=\frac{K_1}{2a_\mu}e^{\sum_{j=1}^nb_jz_j+b_{n+1}}z_\mu+\frac{K_2}{2a_\mu}e^{-\sum_{j=1}^nb_jz_j-b_{n+1}}z_\mu+\beta_\mu,\eeas 
where $\beta_j\in\mathbb{C}$ for $1\leq j\leq n$. Hence the solution is $\Phi(\beta_1,\beta_2,\ldots, \beta_n )=0$. For simplicity, we suppose
\be\label{w13} f_1(z)=\frac{z_\mu}{2a_\mu}\left(K_1e^{\sum_{j=1}^nb_jz_j+A}+K_2e^{-\sum_{j=1}^nb_jz_j-A}\right)+h_1\left(y\right),\ee
where $h_1\left(y\right)$ is a finite order entire function satisfying $\sum_{j=1}^n a_j\frac{\pa h_1(y)}{\pa z_j}\equiv 0$.
Using similar argument as above, we can derive from the third equation in (\ref{teq1}) that
\be\label{w14} f_2(z)=\frac{z_\mu}{2a_\mu}\left(K_3e^{-\sum_{j=1}^nb_jz_j+B}+K_4e^{\sum_{j=1}^nb_jz_j-B}\right)+h_2(y),\ee
where $h_2(y)$ is a finite order entire function  satisfying $\sum_{j=1}^n a_j\frac{\pa h_2(y)}{\pa z_j}\equiv 0$.
By utilizing the equations (\ref{w13}) and (\ref{w14}) in the second and fourth equations of (\ref{teq1}), we get respectively
\beas&& a_{n+1}h_1(y)+a_{n+2}h_2(y+s)
\equiv\gamma_1(n)K_1e^{\sum_{j=1}^nb_jz_j+A}+\gamma_2(n)K_2e^{-\sum_{j=1}^nb_jz_j-A}\\\text{and}
&& a_{n+1}h_2(y)+a_{n+2}h_1(y+s)
\equiv\gamma_1(n)K_3e^{-\sum_{j=1}^nb_jz_j+B}+\gamma_2(n)K_4e^{\sum_{j=1}^nb_jz_j-B},\eeas
where $\gamma_1(n)$ and $\gamma_2(n)$ are given in (\ref{K2}).\end{proof}
\begin{proof}[\bf{Proof of Theorem \ref{T23}}]
Let
\beas\begin{array}{ll}
\frac{K_1}{ia_{n+2}K_4}\Gamma_1(z)e^{P(z)+Q(z+c)}\equiv 1\;\;\text{and}\;\;
\frac{K_4}{ia_{n+2}K_2}\Gamma_2(z)e^{-Q(z)+P(z+c)}\equiv 1.\end{array}\eeas
By using similar argument to those of  \textrm{Theorem \ref{T13}}, we get a contradiction.\\
Let
\beas\begin{array}{ll}
\frac{K_2}{ia_{n+2}K_4}\Gamma_1(z)e^{-P(z)+Q(z+c)}\equiv1\;\;\text{and}\;\;
\frac{K_3}{ia_{n+2}K_2}\Gamma_2(z)e^{Q(z)+P(z+c)}\equiv 1.\end{array}\eeas
By using similar argument to those of  \textrm{Theorem \ref{T13}}, we get a contradiction.\\
Let
\bea\label{w4}\begin{array}{ll}
\frac{K_2}{ia_{n+2}K_4}\Gamma_1(z)e^{-P(z)+Q(z+c)}\equiv1\;\;\text{and}\;\;
\frac{K_4}{ia_{n+2}K_2}\Gamma_1(z)e^{-Q(z)+P(z+c)}\equiv 1.\end{array}\eea
Using (\ref{w4}), from (\ref{mnb4}) and (\ref{mnb5}), we have
\bea\label{w41}\begin{array}{ll}
 \frac{K_1}{ia_{n+2}K_3}\Gamma_1(z)e^{P(z)-Q(z+c)}\equiv 1\;\;\text{and}\;\;
\frac{K_3}{ia_{n+2}K_1}\Gamma_2(z)e^{Q(z)-P(z+c)}\equiv 1.\end{array}\eea
By similar argument as of \textrm{Theorem \ref{T22}}, we have
$P(z)=\sum_{j=1}^nb_jz_j+A$ and $Q(z)=\sum_{j=1}^nb_jz_j+B$, where $b_j,A,B\in\mathbb{C}$ for $1\leq j\leq n$.
From (\ref{w4}) and (\ref{w41}), we have
\be\label{w42}\left\{\begin{array}{ll}
\frac{K_2}{ia_{n+2}K_4}\left(\sum_{j=1}^n a_j b_j-ia_{n+1}\right)e^{\sum_{j=1}^nb_jc_j-A+B}\equiv 1,\\[1mm]
\frac{K_4}{ia_{n+2}K_2}\left(\sum_{j=1}^n a_j b_j-ia_{n+1}\right)e^{\sum_{j=1}^nb_jc_j+A-B}\equiv 1,\\[1mm]
 \frac{K_1}{ia_{n+2}K_3}\left(\sum_{j=1}^n a_j b_j-ia_{n+1}\right)e^{-\sum_{j=1}^nb_jc_j+A-B}\equiv 1,\\[1mm]
\frac{K_3}{ia_{n+2}K_1}\left(\sum_{j=1}^n a_j b_j-ia_{n+1}\right)e^{-\sum_{j=1}^nb_jc_j-A+B}\equiv 1.\end{array}\right.\ee
From (\ref{w42}), we have $\left(\sum_{j=1}^n a_j b_j-ia_{n+1}\right)^2=-a_{n+2}^2$, {\it i.e.,} $\sum_{j=1}^n a_j b_j=i(a_{n+1}\pm a_{n+2})$.\\[1.5mm]
Let $\sum_{j=1}^n a_j b_j=i(a_{n+1}-(-1)^\nu a_{n+2})$ $(\nu=1,2)$. Now two different cases are possible $a_{n+1}\not=(-1)^\nu a_{n+2}$ and $a_{n+1}=(-1)^\nu a_{n+2}$. The second case is consider in the proof of the \textrm{Theorem \ref{T24}}.\\
If $a_{n+1}\not=(-1)^\nu a_{n+2}$, then by similar argument as of \textrm{Theorem \ref{T22}}, we get
\bea\label{w43}\begin{array}{ll}
f_1(z)=\frac{K_1e^{\sum_{j=1}^nb_jz_j+A}-K_2e^{-\sum_{j=1}^nb_jz_j-A}}{2\sum_{j=1}^n a_j b_j}+h_3(y),\\
f_2(z)=\frac{K_3e^{\sum_{j=1}^nb_jz_j+B}-K_4e^{-\sum_{j=1}^nb_jz_j-B}}{2\sum_{j=1}^n a_j b_j}+h_4(y),\end{array}\eea
where $h_k(y)$ $(3\leq k\leq 4)$ are finite order entire functions of $a_1z_2-a_2z_1,\ldots,a_1z_n-a_nz_1$ satisfying $\sum_{j=1}^n a_j\frac{\pa h_k(y)}{\pa z_j}\equiv 0$ and
\bea\label{w44}\left\{\begin{array}{lll}
e^{2\sum_{j=1}^nb_jc_j}\equiv 1, e^{2A-2B}\equiv K_2K_3/K_1K_4, e^{\sum_{j=1}^nb_jc_j-A+B}\equiv (-1)^{\nu+1}K_4/K_2,\\[1.5mm]
e^{\sum_{j=1}^nb_jc_j+A-B}\equiv (-1)^{\nu+1}K_2/K_4, e^{-\sum_{j=1}^nb_jc_j+A-B}\equiv (-1)^{\nu+1}K_3/K_1,\\[1.5mm]
 e^{-\sum_{j=1}^nb_jc_j-A+B}\equiv (-1)^{\nu+1}K_1/K_3\quad\text{where}\quad \nu=1,2.\end{array}\right.\eea
Using (\ref{w43}), (\ref{w44}) and $\sum_{j=1}^n a_jb_j=i(a_{n+1}-(-1)^\nu a_{n+2})$, we deduce from the second and fourth equations of (\ref{teq1}) that
$ a_{n+1}h_3(y)+a_{n+2}h_4(y+s)\equiv 0$ and $a_{n+1}h_4(y)+a_{n+2}h_3(y+s)\equiv 0$.\end{proof}
\begin{proof}[\bf{Proof of Theorem \ref{T24}}]
If $a_{n+1}=(-1)^\nu a_{n+2}$ $(\nu=1,2)$, then $\sum_{j=1}^n a_jb_j=0$. By similar argument to those of \textrm{Theorem \ref{T22}}, we can get the conclusions. So we are omitting the details.
This completes the proof.
\end{proof}
\section{Declarations}
\noindent{\bf Acknowledgment:} The work of the first author is supported by University Grants Commission (IN) fellowship. \\
{\bf Author's contributions:} All authors have equal contribution to complete the manuscript. All of them read and approved the final manuscript.\\
{\bf Conflict of Interest:} Authors declare that they have no conflict of interest.\\
{\bf Availability of data and materials:} Not applicable.\\

\end{document}